\documentclass{amsart}
\usepackage{epsfig}
\usepackage{amsmath}
\usepackage{amsthm}
\usepackage{latexsym}
\usepackage{graphicx}
\usepackage{amsfonts}
\usepackage{amssymb}
\newtheorem{theorem}{Theorem}[section]
\newtheorem{lemma}[theorem]{Lemma}

\newtheorem{proposition}[theorem]{Proposition}

\begin{document}

\title[The maximal number of exceptional Dehn surgeries]
{The maximal number\\ of exceptional Dehn surgeries}

\author{Marc Lackenby}
\thanks{First author supported by an EPSRC Advanced Research Fellowship.}
\author{Robert Meyerhoff}
\thanks{Second author partially supported by NSF grants DMS-0553787 and DMS-0204311.}

\address{Mathematical Institute\\University of Oxford\\Oxford, UK}
\address{Department of Mathematics\\Boston College\\Chestnut Hill, MA}

\maketitle

\section{Introduction}

Thurston's hyperbolic Dehn surgery theorem is one of the most important results
in 3-manifold theory, and it has stimulated an enormous amount of research. If 
$M$ is a compact orientable hyperbolic 3-manifold with boundary a single torus,
then the theorem asserts that, for all but finitely many slopes $s$ on $\partial M$, the
manifold $M(s)$ obtained by Dehn filling along $s$ also admits a
hyperbolic structure. The slopes $s$ where $M(s)$ is not hyperbolic
are known as {\sl exceptional}. A major open question has been: what is
the maximal number of exceptional slopes on such a manifold $M$? 
When $M$ is the exterior of the figure-eight knot, the number of exceptional
slopes is 10, and this was conjectured by Gordon in \cite{Gor1} to be an
upper bound that holds for all $M$. In this paper, we prove this
conjecture.

\begin{theorem}
Let $M$ be a compact orientable
3-manifold with boundary a torus, and with interior admitting
a complete finite-volume hyperbolic structure. Then the number
of exceptional slopes on $\partial M$ is at most $10$.
\label{thm:main1}\end{theorem}

Although it is not our approach,
much of the progress on this problem has been achieved by
bounding the intersection number $\Delta(s_1,s_2)$ between
exceptional slopes $s_1$ and $s_2$. It was conjectured by Gordon \cite{Gor1}
that $\Delta(s_1, s_2)$ is always at most $8$, a bound which is attained for
the exteriors of the figure-eight knot and the figure-eight knot sister. We also prove this conjecture.

\begin{theorem}
Let $M$ be a compact orientable
3-manifold with boundary a torus, and with interior admitting
a complete finite-volume hyperbolic structure. If $s_1$ and
$s_2$ are exceptional slopes on $\partial M$, then their
intersection number $\Delta(s_1, s_2)$ is at most $8$.
\label{thm:main2}\end{theorem}

Recently and using different methods, Agol \cite{Ago2} has shown that, for all but finitely
many 3-manifolds $M$ as in Theorem 1.1, the intersection number between exceptional
slopes is at most $5$ and the number of exceptional slopes is at most $8$.
Moreover, there is an algorithm to compute the list of excluded manifolds.
However, the algorithm is far from practical, and so there seems, at
present, to be no way of using Agol's theorem to prove Gordon's
conjectures.

Theorems \ref{thm:main1} and \ref{thm:main2} are proved using a combination of new geometric
techniques and a rigorous computer-assisted calculation. In particular,
an improved version of the 6-theorem of Lackenby  \cite{Lac1}
and Agol \cite{Ago1} is established, and extensive use of the Mom technology of
Gabai, Meyerhoff and Milley (\cite{GMM2}, \cite{GMM3}) is required.

\vfill\eject
In Mom theory, the following list of 2-cusped and 3-cusped hyperbolic 3-manifolds
plays a central role, where the notation is that of the
census \cite{CHW1} of hyperbolic 3-manifolds:

\begin{figure}[htp]
\begin{center}
\begin{tabular}{|c|c|c|c|c|c|c|c|c|}
\hline\hline
$m412$ & $s596$ & $s647$ & $s774$ & $s776$ & $s780$ & $s785$ & $s898$ & $s959$ \\
\hline\hline
\end{tabular}
\end{center}
\label{fig:mommflds}
\caption{}
\end{figure}

The proof of Theorems \ref{thm:main1} and \ref{thm:main2}
immediately divides into two cases: either $M$ is obtained by Dehn filling
one of the manifolds in this list, or it is not. In the former case, 
a straightforward analysis using the computer program Snap \cite{Goo1} leads to a proof
of the theorems (see Section 8). 
The other case is when $M$ is not obtained by Dehn filling
one of the manifolds in this list. We consider the inverse image in ${\Bbb H}^3$ of a maximal horoball
neighbourhood of the cusp of $M$, which is a collection of horoballs.
Following \cite{GMM2} and \cite{GMM3}, we extract three real-valued parameters from
this arrangement, which are denoted $e_2$, $e_3$ and $e_4$.
(More details can be found in Section 3.) We also define three
other real-valued parameters, $m$, $t$ and $h$, which encode
the shape and size of the cusp torus. These 6 parameters then
define a parameter space.
We show that outside an explicit compact region of this
parameter space, the theorems hold. We then examine this compact region
using a rigorous computer analysis. We divide the region into small pieces,
and show that, in each piece, the theorems hold. This requires two approaches. 
We develop new geometric tools which can be
used to deduce that, in certain regions, either a contradiction is reached
or $M$ contains a `torus-friendly geometric Mom-2 or Mom-3' and hence
is obtained by Dehn filling one of the manifolds in Figure \ref{fig:mommflds}. 
If we cannot exclude a region, then we need to bound the number of exceptional
slopes in the boundary torus and their intersection numbers. 
Prior to the present paper, this amounted to using the 6-theorem and checking whether or not the length of the 
slope was at most 6.  However, for our purposes, this is not strong enough and so
we prove an enhanced version of the 6-theorem, which depends on the parameter $e_2$.

The plan of the paper is as follows. In Section 2, we give a
historical survey of the problem. In Section 3, we recall
the terminology and techniques of Mom structures. In Section 4, we
prove the extension of the 6-theorem. In Section 5, we derive
lower bounds on the area of the cusp torus.
In Section 6, we give the geometric arguments that underpin the
parameter space analysis. In Section 7, we give details of the rigorous computation. 
In Section 8, we examine the manifolds that are obtained
by Dehn filling one of the manifolds from Figure \ref{fig:mommflds}, and prove the theorems in this case.

\section{Historical survey}

There have been at least four parallel
approaches to the Dehn surgery problem. The most popular method, with the most extensive literature,
has been topological. Here, one seeks not to establish that the filled-in manifold
$M(s)$ has a hyperbolic structure, but rather that it is `hyperbolike'. The precise
definition of this term varies according to the context, but the gist is that
a compact orientable 3-manifold is hyperbolike if it has topological properties that are equivalent to
the existence of a hyperbolic structure, assuming the geometrisation conjecture.
Thus, now that Perelman's proof of this conjecture (\cite{Per1}, \cite{Per2}, \cite{Per3}) 
is accepted as correct \cite{MT1}, a
compact orientable 3-manifold is hyperbolike if and only if it is hyperbolic. In the topological
approach, $M(s)$ is hyperbolike if it is irreducible, atoroidal and not a
Seifert fibre space. One considers slopes $s_1$ and $s_2$ where neither
$M(s_1)$ nor $M(s_2)$ is hyperbolike and one seeks to bound their intersection
number $\Delta(s_1, s_2)$. For example, if $M(s_1)$ and $M(s_2)$ are both reducible,
then they contain essential 2-spheres, which restrict to essential planar
surfaces in $M$. One considers the intersection pattern of these surfaces and,
after some subtle and beautiful combinatorics, very accurate bounds on
$\Delta(s_1, s_2)$ can be achieved. In fact, in this case, Gordon
and Luecke \cite{GL1} proved that $\Delta(s_1, s_2) \leq 1$, and so there are at most
3 slopes $s$ on $M$ for which $M(s)$ is reducible. In addition, the ${\rm SL}(2, {\Bbb C})$
representation variety has been used extensively in this approach. Many mathematicians have been
part of this program, including Boyer, Culler, Gordon, Luecke, Scharlemann, Shalen,
Wu and Zhang (see \cite{Boy1} for a survey).
However, the case where $M(s_1)$ or $M(s_2)$ is an atoroidal Seifert fibre space
has proved to be problematic. When it has finite fundamental group, the methods
of Boyer and Zhang \cite{BZ1} have been very successful, but when it does not, the situation is
harder to handle topologically, and little progress has been made.

A second approach has been via foliations and laminations. Here, the aim
is to find an essential foliation or lamination on $M(s)$. By Novikov's theorem \cite{Nov1},
this implies that $M(s)$ is irreducible or $S^2 \times S^1$ (and one can typically
rule out the latter case). Moreover, if the essential lamination is genuine, 
then $M(s)$ is hyperbolike, by Gabai and Kazez's theorem \cite{GK1}. According
to a result of Gabai \cite{Mos1}, there is a slope $\lambda$ such that
$M(s)$ has a genuine essential lamination provided $\Delta(s,\lambda) \geq 3$. But unfortunately,
the number of excluded slopes, where $\Delta(s,\lambda) \leq 2$, is not finite. Nevertheless, the foliation
and lamination approach has been extremely successful, notably with Gabai's
proof \cite{Gab1} of the Property R conjecture.

A third approach, due to Hodgson and Kerckhoff \cite{HK1}, aims to prove that $M(s)$
is hyperbolic without appealing to the geometrisation conjecture. Since their
methods predate Perelman's work, this was the first approach to yield a
universal upper bound (60) on the number of exceptional slopes, independent
of the manifold $M$.

The fourth approach has also been geometric, but with the aim of deducing a
weaker conclusion on $M(s)$ than hyperbolicity. The first result in this
direction was the Gromov-Thurston $2\pi$-theorem \cite{BH1} which established a universal upper bound on the
number of slopes $s$ for which $M(s)$ does not admit a Riemannian
metric of negative curvature.  This upper bound, due to Thurston,
was 48. However, there were successive improvements to this bound.
Bleiler and Hodgson reduced it to 24 in \cite{BH1}. By applying work of
Cao and Meyerhoff \cite{CM1}, this could be reduced to 14. 

A continuation of this approach was the $6$-theorem of Lackenby \cite{Lac1} and
Agol \cite{Ago1}, which reduced the number of excluded slopes to $12$. 
However, outside this set of excluded slopes,
$M(s)$ is shown only to be irreducible, atoroidal, and not Seifert
fibered, and to have infinite, word-hyperbolic fundamental group.
But with the solution of the geometrisation conjecture by Perelman, 
any compact orientable 3-manifold satisfying these conditions must also
admit a hyperbolic structure. Thus, as a consequence of Perelman's work, the number of 
exceptional slopes on $M$ is established by this approach to be at most 12.

The $2 \pi$-theorem and $6$-theorem both are
phrased in terms of the length of a slope $s$. Here, one considers
the unique maximal horoball neighbourhood of the cusp of $M$, and
its immersed boundary torus, which we term the {\sl cusp torus}. One defines the {\sl length} of $s$ to be the 
length of shortest curve on this torus with slope $s$. If the length of $s$ 
is more than $2 \pi$, then $M(s)$ admits a negatively curved Riemannian metric.
If the length is more than $6$, then $M(s)$ is irreducible, atoroidal, and not Seifert
fibered, and has infinite, word-hyperbolic fundamental group.
Thurston gave an upper bound of $48$ on the
number of slopes with length at most $2 \pi$. 

The argument that leads to the bound of 48 is instructive. Using elementary means, Thurston showed
that the length of each slope is at least 1. Thus, it is
clear that, for any fixed constant $L$, there is a uniform upper
bound on the number of slopes with length at most $L$. To actually
find this upper bound, one can argue as follows. The universal cover
of the cusp torus is a Euclidean plane, and the inverse image of a basepoint
is a lattice in this plane. The fact that each slope has length at least 1 forces
every pair of lattice points to be at least distance 1 apart.
Thus, the minimal possible co-area of such a lattice is $\sqrt 3/2$, which is
achieved by the hexagonal lattice. The area $A$ of the cusp torus is therefore at
least $\sqrt 3/2$. An elementary argument (see the proof of Theorem 8.1 in 
\cite{Ago1} for example) gives that
the intersection number of two slopes with length $\ell_1$ and $\ell_2$
is at most $\ell_1 \ell_2 / A$. When $\ell_1, \ell_2 \leq 2 \pi$
and $A \geq \sqrt 3/2$, the intersection number is therefore at most 45.
Lemma 8.2 in \cite{Ago1} states that, if $E$ is any collection of slopes on a torus,
where any two slopes in $E$ have intersection number at most $\Delta$, then
$|E| \leq p+1$, where $p$ is any prime more than $\Delta$. (Note that this bound is
not always sharp. For example, if $\Delta = 7$, then $|E|$ can be shown to be
at most $10$.) Setting $p=47$ gives the upper bound.

Thus, the formula $\ell_1 \ell_2 / A$ is central to the fourth approach
to the Dehn surgery problem. Increasingly good
upper bounds on $\ell_1 \ell_2 / A$, when $\ell_1$ and $\ell_2$ are
the lengths of exceptional slopes, have been found. The $6$-theorem reduced the critical
value of $\ell_1$ and $\ell_2$ from $2 \pi$ to $6$. But, the most significant
progress has been achieved by finding improved lower bounds on $A$, 
the area of the cusp torus. Adams \cite{Ada1} increased the lower
bound on $A$ to $\sqrt 3$, which gave the upper
bound of 24 exceptional slopes. Then Cao and Meyerhoff \cite{CM1} gave a lower
bound of 3.35 for $A$, which led to the upper bound of 12
exceptional slopes. Recently, the work of Gabai, Meyerhoff and Milley (\cite{GMM2}, \cite{GMM3})
improves the lower bound for $A$ to 3.7, provided that $M$ is not obtained
by Dehn filling one of the manifolds in Figure \ref{fig:mommflds}. This result is
not explicitly stated in their paper, but it follows from their methods.
This leads to an upper bound of $9$ on the intersection number of any
two exceptional slopes, but unfortunately, it does not decrease the bound
on the number of exceptional slopes below 12. In fact, the following table, taken from \cite{Boy1}, gives
the maximal size of a set of slopes, such that any two slopes in this set
have intersection number at most $\Delta$.

\vskip 6pt
\begin{center}
\begin{tabular}{|c|c|c|c|c|c|c|c|c|c|c|c|}
\hline\hline
Maximal intersection number $\Delta$ 
& 0 & 1 & 2 & 3 & 4 & 5 & 6 & 7 & 8 & 9 & 10 \\
\hline
Maximal number of slopes 
& 1 & 3 & 4 & 6 & 6 & 8 & 8 & 10 & 12 & 12 & 12 \\
\hline\hline
\end{tabular}
\end{center}
\vskip 6pt

So, reducing the upper bound on the number of exceptional slopes 
from $12$ to $10$ is not as straightforward as it first appears. If one were to
prove the bound of $10$ using intersection numbers, then one would have to
improve the known bounds on intersections numbers from $9$ to $7$, which
is quite a significant reduction.

In addition, there is an example of Agol that illustrates some difficulties.
If one wants to develop the fourth approach to the problem,
it seems that either one must improve the $6$-theorem yet further,
by reducing the critical slope length below $6$, or one must show that $M$ always has
at most 10 slopes with length at most 6. But neither of these
steps is possible. Agol \cite{Ago1} examined the case where $M$ is the exterior of
figure-eight knot sister. This has two key properties: it has
12 slopes with length at most 6, {\sl and} it has exceptional slopes
with length precisely 6. Fortunately, this manifold
is obtained by Dehn filling $s776$, which appears in Figure \ref{fig:mommflds}. 
Thus, by work of Gabai, Milley and Meyerhoff (\cite{GMM2}, \cite{GMM3}) 
it falls into a family of well-understood 
exceptions. In fact, it is by melding the Mom technology with some new geometric
results, specifically tailored to the Dehn surgery problem, that we are
able to prove Theorems \ref{thm:main1} and \ref{thm:main2}.

\section{Mom terminology and results}

The Mom technology of Gabai, Meyerhoff and Milley plays a key role in this
paper. Several of its main ideas played an implicit part in the
work of Cao and Meyerhoff \cite{CM1}, and we have seen that this was one
of the pieces of the jigsaw that gave the upper bound of
12 exceptional slopes. In this section, we recall the Mom
terminology and results.

Let $M$ be a compact orientable 3-manifold with boundary a torus,
and with interior admitting a complete finite-volume hyperbolic structure.
The universal cover of $M - \partial M$ is hyperbolic 3-space ${\Bbb H}^3$, for which we use the
upper half-space model. The inverse image in ${\Bbb H}^3$ of a
maximal horoball neighbourhood of the cusp is a union of horoballs
$\{ B_i \}$. We may arrange that one of these horoballs, denoted $B_\infty$, is 
$\{ (x,y,z) : z \geq 1 \}$. Given two horoballs $B_i$ and
$B_j$, neither equal to $B_\infty$, we say that $B_i$ and $B_j$ are
in the same {\sl orthoclass} if either $B_i$ and $B_j$ differ
by a covering transformation that preserves $B_\infty$ or
there exists some covering transformation $g$ such that
$g(B_i) = B_\infty$ and $g(B_\infty) = B_j$. We denote the
orthoclasses by ${\mathcal O}(1)$, ${\mathcal O}(2)$, and so
forth. For any $B \in {\mathcal O}(n)$, we call $d(B, B_\infty)$
the {\sl orthodistance} of $B$ and denote it $o(n)$. We order
the orthoclasses in such a way that the corresponding
orthodistances are non-decreasing. Since the horoball neighbourhood of the cusps is
maximal, there is one $B_i$ that touches $B_\infty$, and
hence $o(1) = 0$. It is convenient also
to define the quantity $e_n = \exp(o(n)/2)$. Thus
$e_n$ is a non-decreasing sequence starting at $e_1 = 1$.
In the upper half-space model, the Euclidean diameter
of $B \in {\mathcal O}(n)$ is $e_n^{-2}$. The {\sl orthocentre} of
a horoball $B_i$ other than $B_\infty$ is the closest point on
$B_\infty$ to $B_i$. We say that two horoballs $B_i$ and $B_j$ are
in the ${\mathcal O}(n)$ {\sl orthopair class} if there is a covering 
transformation $g$ such that $g(B_i) = B_\infty$ and $g(B_j) \in {\mathcal O}(n)$.
We then say that $B_i$ and $B_j$ are {\sl $o(n)$-separated}.

A {\sl $(p,q,r)$-triple} is a triple of horoballs $\{ B_1, B_2, B_3 \}$
such that $B_1, B_2$ are in the ${\mathcal O}(r)$ orthopair class,
$B_2, B_3$ are in the ${\mathcal O}(p)$ orthopair class and $B_1, B_3$ are in
the ${\mathcal O}(q)$ orthopair class, possibly after re-ordering. A {\sl geometric Mom-$n$ structure}
is a collection of $n$ triples of type $(p_1, q_1, r_1), \dots,
(p_n, q_n, r_n)$, no two of which are equivalent under the
action of $\pi_1(M)$, and such that the indices $p_i$, $q_i$
and $r_i$ all come from the same $n$-element subset of
${\Bbb Z}_+$. We say that the Mom-$n$ structure {\sl involves} ${\mathcal O}(k)$ if 
$k \in \{ p_1, q_1, r_1, \dots, p_n, q_n, r_n \}$. A geometric
Mom-$n$ is {\sl torus-friendly} if $n=2$ or
if $n=3$ and the Mom-3 does not possess exactly two triples
of type $(p,q,r)$ for any set of distinct positive
indices $p$, $q$ and $r$.

The first key lemma, which appears as Corollary 3.3 in \cite{GMM3}, is as follows.

\begin{lemma}
For any integer $n$, there are no $(n,n,n)$-triples.
\label{lem:nonnn}\end{lemma}

In particular, there are no $(1,1,1)$-triples. Hence, the Euclidean
distance in $\partial B_\infty$ between the orthocentres of
two ${\mathcal O}(1)$-horoballs is at least $e_2$. More generally,
one can compute the distance between orthocentres using the following
elementary lemma (Lemma 3.4 in \cite{GMM3}).

\begin{lemma}
Suppose that $B_1 \in {\mathcal O}(q)$,
$B_2 \in {\mathcal O}(r)$, and that $B_1$ and $B_2$ are $o(p)$-separated.
Then the Euclidean distance between the orthocentres of $B_1$ and $B_2$
is $e_p/(e_q e_r)$.
\label{lem:orthocentreseparation}
\end{lemma}

The following results will be crucial. These are contained in \cite{GMM3} but not explicitly
stated there. We include here only a brief outline of their proof. 

\begin{theorem}
If $M$ contains a  geometric Mom-2 involving ${\mathcal O}(1)$ and ${\mathcal O}(n)$, 
then either $e_n \geq 1.5152$ or $M$ is obtained by Dehn filling
one of the manifolds $m125$, $m129$ or $m203$.
\label{thm:mom2mflds}\end{theorem}

Note that the manifolds $m125$, $m129$ or $m203$ are obtained by Dehn filling
$s776$, which appears in Figure \ref{fig:mommflds}.

\begin{theorem}
If $M$ contains a torus-friendly geometric Mom-3 involving ${\mathcal O}(1)$,
${\mathcal O}(2)$ and ${\mathcal O}(3)$, then either $e_3 \geq 1.5152$ or $M$ is obtained by Dehn filling
one of the manifolds in Figure \ref{fig:mommflds}.
\label{thm:mom3mflds}\end{theorem}

\noindent\emph{Proof:} 
Suppose first that $M$ contains a geometric Mom-2 involving ${\mathcal O}(1)$ and ${\mathcal O}(n)$. 
Associated to this, there is  a cell complex $\Delta$, defined in \cite{GMM3}. Since $e_n \leq 1.5152$, the material in
Section 6 of \cite{GMM3} gives that $\Delta$ is embedded. Then Theorem 7.1
and Proposition 8.2 of \cite{GMM3} give that $M$ is obtained by Dehn filling
a `full topological internal Mom-2 structure'. By Theorem 5.1 of \cite{GMM2},
there are only 3 hyperbolic manifolds with such a structure: $m125$, $m129$ and $m203$.

Suppose now that $M$ contains a torus-friendly geometric Mom-3 involving
${\mathcal O}(1)$, ${\mathcal O}(2)$ and ${\mathcal O}(3)$
but not a geometric Mom-2 involving a subset of ${\mathcal O}(1)$, ${\mathcal O}(2)$ and ${\mathcal O}(3)$.
Then, because $e_3 \leq 1.5152$, the associated cell complex is embedded.
Theorem 7.1 and Theorem  8.3 of \cite{GMM3} gives that $M$ is obtained by Dehn filling
a `full topological internal Mom-3 structure'. By Theorem 5.1 of \cite{GMM2},
any hyperbolic manifold admitting such a structure is obtained by Dehn filling one of the
manifolds in Figure \ref{fig:mommflds}. \qed\

\section{An improvement to the 6-theorem}

In this section, we prove the following result, which is a version
of the 6-theorem that depends on the parameter $e_2$. 

\begin{theorem}
Let $M$ be a compact orientable
3-manifold, with boundary a torus and with interior admitting
a complete finite-volume hyperbolic structure. Let $s$
be a slope on $\partial M$ with length at least
$${\pi e_2 \over \arcsin(e_2/2)}$$
if $e_2 \leq \sqrt 2$, and length at least
$${2\pi e_2 \over 2 \arcsin (\sqrt{ 1 - e_2^{-2} }) + e_2^2 - 2 \sqrt{ e_2^2 - 1 } }$$
if $e_2 > \sqrt 2$. Then, $M(s)$ is hyperbolic.
\label{thm:improve6}\end{theorem}

\begin{figure}[htp]
\begin{center}
\includegraphics{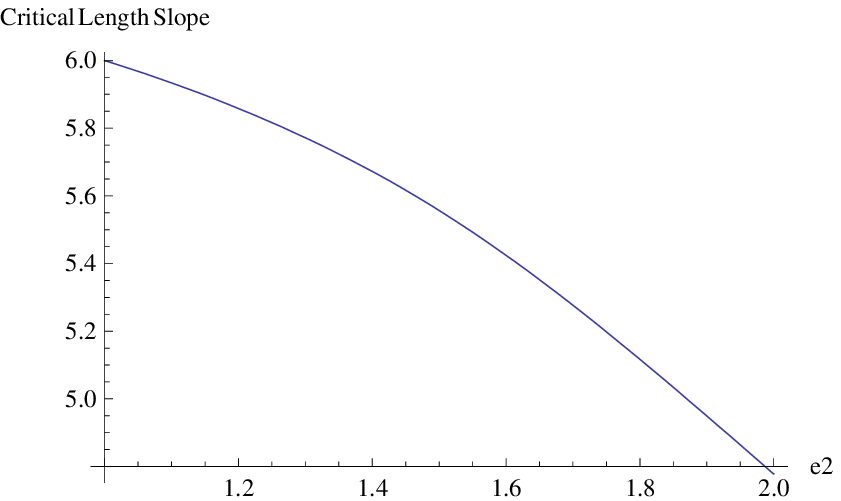}
\end{center}
\caption{}
\label{fig:graph}
\end{figure}

The critical slope length, as a function of $e_2$, is as shown in Figure \ref{fig:graph}.
When $e_2 = 1$, Theorem \ref{thm:improve6}
gives the same critical slope length as the 6-theorem. But as $e_2$ increases,
the critical slope length decreases, tending to zero.

We start by recalling an outline of the proof of the 6-theorem. The proof of
Theorem \ref{thm:improve6} will be a refinement of this. Suppose that $M(s)$ is
not hyperbolic. We wish to show that the length of $s$ 
is at most $6$. Suppose, for simplicity, that $M(s)$ contains an essential sphere.
By choosing this sphere suitably, we may arrange that its intersection $F$ with $M$
is incompressible and boundary-incompressible and that each boundary
component of $F$ has slope $s$. We may then homotope $F$ to a pleated surface.
Its area is then at most $-2 \pi \chi(F) = 2 \pi (|\partial F| - 2)$.
The aim is to show geometrically that each component of $\partial F$ contributes at least
$(\pi /3)$ times the length of $s$ to the area of $F$. Thus, if the
length of each $s$ is more than $6$, then the area of $F$ is at least
$(\pi /3) 6 |\partial F| = 2 \pi |\partial F|$, which is a contradiction.
The intersection of $F$ with the maximal horoball neighbourhood of the cusps is
a collection of copies of $S^1 \times [1,\infty)$ plus possibly some compact components.
Each copy of $S^1 \times [1,\infty)$ is ambient isotopic to a vertical surface lying over
a curve of slope $s$ in the cusp torus. Its area is
therefore at least the length of $s$. Thus, each component of $\partial F$ contributes
at least ${\rm Length}(s)$ to the area of $F$. However, we want to improve this
contribution to $(\pi /3) {\rm Length}(s)$, and to do this, one must consider the
parts of $F$ not lying in the horoball neighbourhood of the cusps. The most convenient way to do this
is to consider the associated geodesic spine $S$, which is defined as follows. 

The inverse image in ${\Bbb H}^3$
of the maximal horoball neighbourhood of the cusp is a collection $\{ B_i \}$ of
horoballs. One of these, $B_\infty$, has been fixed as $\{ (x,y,z) : z \geq 1 \}$ in the upper half-space 
model. Let $\tilde S$ be the set of points in ${\Bbb H}^3$ that do not have a unique
closest point in $\bigcup B_i$. It is invariant under the group of
covering transformations, and its quotient in $M$ is $S$. This is a spine
for $M$, in which each cell is totally geodesic. Thus, $M - S$ is a neighbourhood
of the cusp that is larger than the interior of the maximal horoball neighbourhood.
By considering the area of $F$ in $M - S$, we obtain the improved area
contribution of $(\pi /3) {\rm Length}(s)$. Specifically, the following
result is used, which appears as Lemma 3.3 in \cite{Lac1}.

\begin{lemma}
Let $M$ be a compact orientable 3-manifold, with boundary a torus and with interior
admitting a complete finite-volume hyperbolic structure. Let $S$ be a geodesic spine
arising from a horoball neighbourhood of the cusp of
$M - \partial M$. Let $G$ be a compact orientable
(possibly non-embedded) surface with interior in $M - S$, with boundary in 
$\partial M \cup S$ and with $\partial G \cap \partial M$
representing $\pm k [s] \in H_1(\partial M)$,
where $k \in {\Bbb N}$ and $s$ is some slope.
Then
$${\rm Area}(G -  \partial M) \geq k \, (\pi /3) {\rm Length}(s).$$
\label{lem:wordareabound}\end{lemma}

By applying this to each component of the surface $F - S$, we obtain the 6-theorem, at least
when $M(s)$ is reducible.

If $M(s)$ has finite fundamental group, then the core of the surgery
solid torus has finite order. So, some power of this core curve bounds a
disc in $M(s)$. The restriction of this disc to $M$ is a compact planar
surface $F$, with all but one boundary component having slope $s$.
We apply the above argument to this surface.

To show that $\pi_1(M(s))$ is word hyperbolic, we use Gabai's ubiquity theorem \cite{Gab2},
as stated as Theorem 2.1 in \cite{Lac1}. We consider an arbitrary loop $K$ in $M$ that
is homotopically trivial in $M(s)$. There is therefore a compact planar
surface $F$, with one boundary component mapped to $K$, and the remaining
boundary components mapped to non-zero multiples of the slope $s$.
We may assume that $F$ is homotopically incompressible and homotopically
boundary-incompressible (as defined in \cite{Lac1}).
The aim is to show that $|F \cap \partial M|$ is bounded above by
$c \ {\rm Length}(K)$, where $c$ is a constant depending on $M$ and $s$, but not $K$.
Theorem 2.1 in \cite{Lac1} then implies that $\pi_1(M(s))$ is word hyperbolic.
In order to establish this bound, we again consider the area of $F$, but
the argument is slightly more complicated. The detailed proof appears
in \cite{Lac1}.

Finally, if $M(s)$ is irreducible, and has infinite, word hyperbolic fundamental
group, then it is atoroidal and not Seifert fibred. Thus, the 6-theorem is
established.

Let us now define a function $I(e_2)$, which will turn out to be the improvement factor in
the critical slope length, as compared with the $2 \pi$ theorem. When $e_2 \leq \sqrt 2$,
$$I(e_2) = {2 \arcsin(e_2/2) \over e_2}.$$
When $e_2 > \sqrt 2$, 
$$I(e_2) = {2\arcsin(\sqrt { 1 - (1/e_2^2) }) + e_2^2 - 2 \sqrt { e_2^2 - 1} \over e_2}.$$
In order to prove Theorem \ref{thm:improve6}, we need the following proposition, which
is an improvement on Lemma \ref{lem:wordareabound}.

\begin{proposition}
Let $M$, $S$, $G$ and $k$ be as in Lemma \ref{lem:wordareabound}.
Then the area of $G -  \partial M$ is at least $k \, {\rm Length}(s) I(e_2)$.
\label{pro:newareabound}\end{proposition}

Thus, the critical slope length is improved to $2 \pi / I(e_2)$,
thereby proving Theorem \ref{thm:improve6}

\begin{figure}[htp]
\begin{center}
\includegraphics{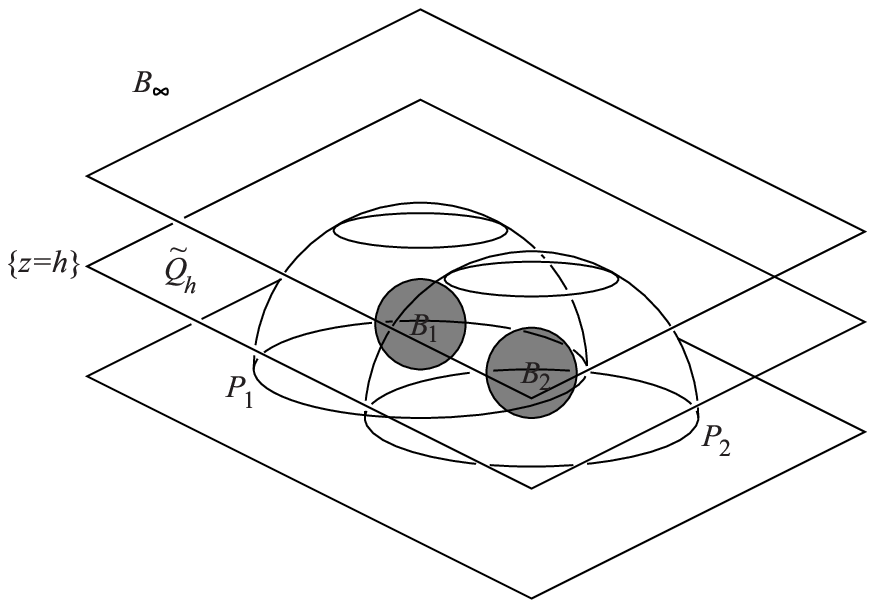}
\end{center}
\caption{}
\label{fig:horoballdefns1}
\end{figure}

\begin{figure}[htp]
\begin{center}
\includegraphics{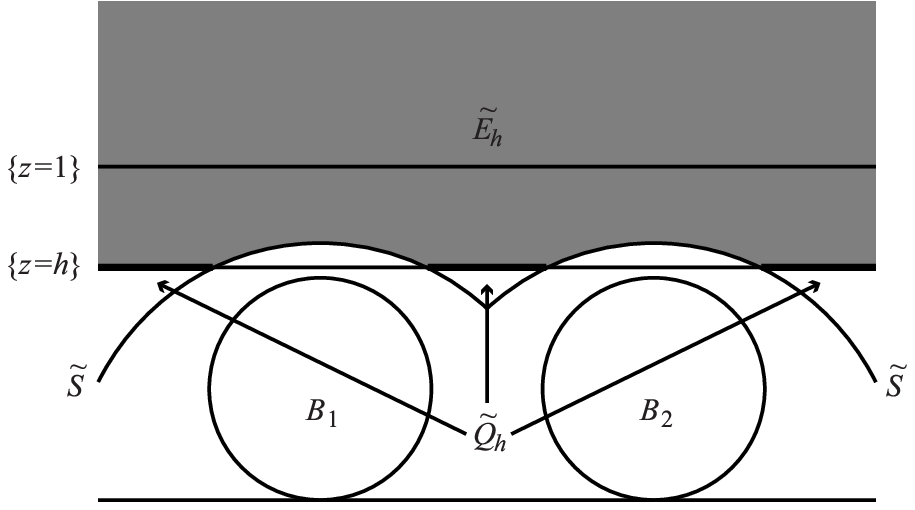}
\end{center}
\caption{}
\label{fig:horoballdefns2}
\end{figure}

We now briefly explain the proof of Lemma \ref{lem:wordareabound}, as its main ideas will
be used in the proof of Proposition \ref{pro:newareabound}. We start by introducing some terminology.
Let $M$, $S$, $G$ and $k$ be as in Lemma \ref{lem:wordareabound}. 
Recall that $\tilde S$ is the inverse image of $S$ in ${\Bbb H}^3$.
Let $\tilde E$ be the closure of the component of
${\Bbb H}^3 - \tilde S$ that contains $B_\infty$. For each horoball $B_i$
of $\{ B_i \}$ other than $B_\infty$, let $P_i$
be the totally geodesic plane equidistant between
$B_\infty$ and $B_i$. Let $h$ be a positive
real number. Let $\tilde Q_h$ be the set of points in
$\{ z = h \}$ that lie above $\tilde S$. (See Figure \ref{fig:horoballdefns1}.) 
Let $Q_h$ be the quotient of $\tilde Q_h$ by the
stabiliser of $B_\infty$. Then $Q_h$ is an
embedded surface in $M - S$, and $\bigcup_{h \in (0, \infty)} Q_h = M - S$.
Let $\tilde E_h$ be the set of points in ${\Bbb H}^3$
that lie above $\tilde S \cup \{ z = h \}$,
and let $E_h$ be the quotient of $\tilde E_h$ by
the stabiliser of $B_\infty$. Thus,
$E_h$ is the set of points in $M - S$ between
$Q_h$ and the cusp. (See Figure \ref{fig:horoballdefns2}.)

\begin{figure}[htp]
\begin{center}
\includegraphics{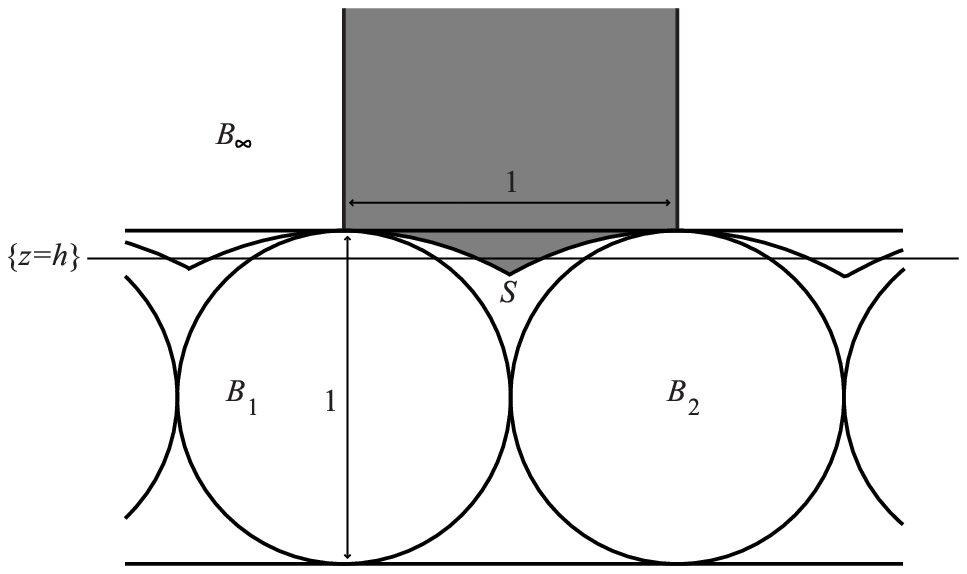}
\end{center}
\caption{}
\label{fig:fig8}
\end{figure}

We may assume, by performing an arbitrarily small perturbation of the surface
$G$, that, for all but finitely many values of $h$, $G \cap Q_h$
is a finite collection of immersed arcs.
The area of $G$ is at least
$$\int_0^\infty {{\rm Length}(G \cap Q_h) \over h} \, dh.$$
The proof of Lemma \ref{lem:wordareabound} proceeded by examining ${\rm Length}(G \cap Q_h) / (k \ {\rm Length}(s))$.
It was shown that this ratio is minimised by a certain surface in
the figure-eight knot complement. In this case, the orthocentres of the ${\mathcal O}(1)$ horoballs
form a hexagonal lattice. Taking a vertical slice
through these, we see an arrangement of horoballs and spine as in Figure \ref{fig:fig8}.
In the proof of Lemma \ref{lem:wordareabound}, it was shown that the ratio 
${\rm Length}(G \cap Q_h) / (k\ {\rm Length}(s))$ is minimised by the intersection of $\{ z = h \}$ 
with the shaded surface. Hence, integrating with respect to $h$, we
get that ${\rm Area}(G) / (k \ {\rm Length}(s))$ is at least that
of the shaded surface. But this surface has area $\pi /3$ and contributes
1 to slope length. This is why the ratio $\pi/3$ appears in Lemma \ref{lem:wordareabound}.

\begin{figure}[htp]
\begin{center}
\includegraphics{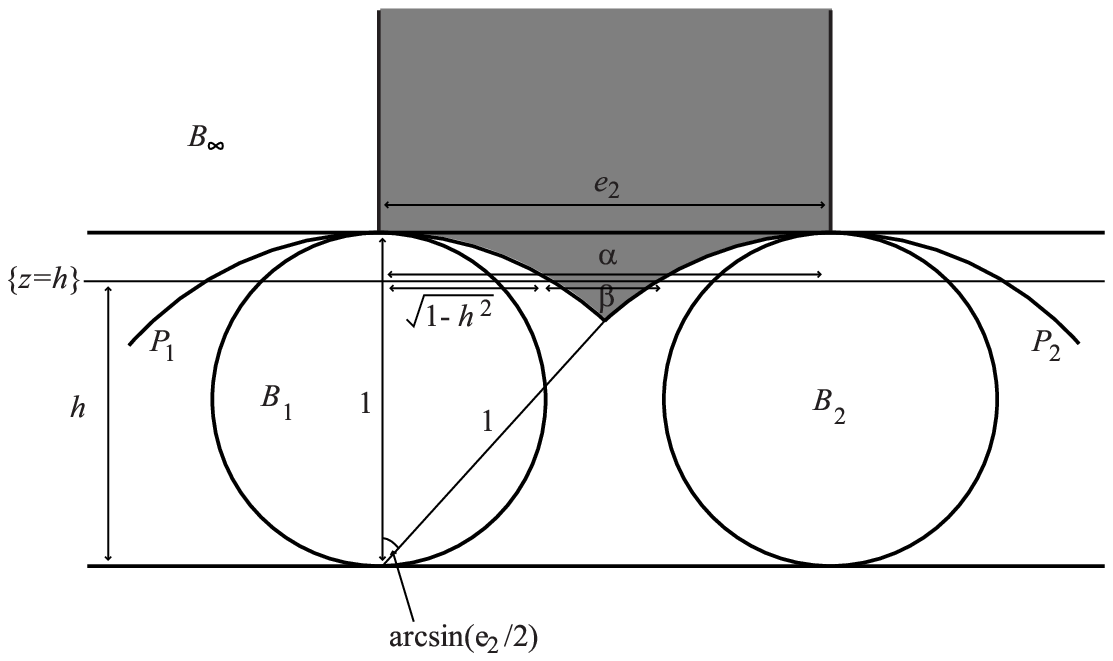}
\end{center}
\caption{}
\label{fig:rellen1}
\end{figure}

We now explain the proof of Proposition \ref{pro:newareabound}. When $e_2$ is bigger than 1, the
${\mathcal O}(1)$ horoballs cannot be tangent. In fact, their orthocentres
are at least $e_2$ apart in the Euclidean metric on $\partial B_\infty$.
Hence, one should be considering an arrangement as in Figure \ref{fig:rellen1}.
The area of the shaded region is $2 \arcsin(e_2/2)$ and its contribution to slope
length is $e_2$. The ratio of these quantities is $I(e_2)$, when $e_2 \leq \sqrt 2$.

\begin{figure}[htp]
\begin{center}
\includegraphics{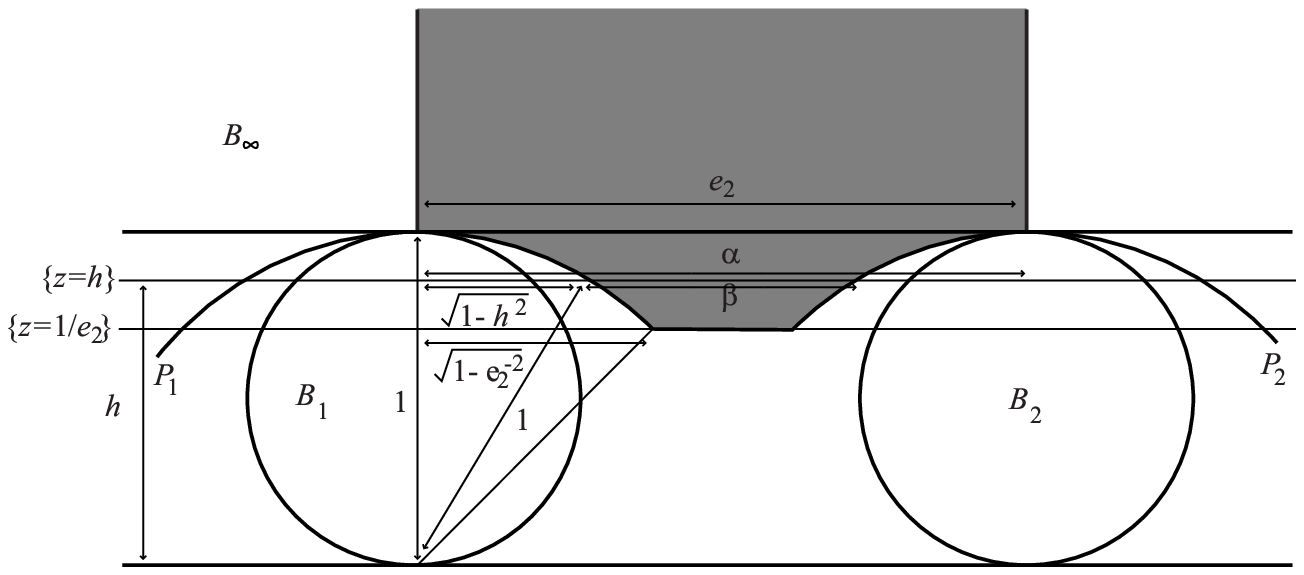}
\end{center}
\caption{}
\label{fig:e2morethanroot2}
\end{figure}

When $e_2 > \sqrt 2$, we again must consider two ${\mathcal O}(1)$ horoballs, but due to the
intervention of other smaller horoballs, it turns out that we must restrict attention to
points above $\{ z = 1/e_2 \}$. Thus, we must consider a configuration as in
Figure \ref{fig:e2morethanroot2}. Here, the shaded region has area
$$2\arcsin(\sqrt { 1 - (1/e_2^2) }) + e_2^2 - 2 \sqrt { e_2^2 - 1}$$
and it contributes $e_2$ to slope length. Thus, again, the ratio of these quantities is
$I(e_2)$. To prove Proposition \ref{pro:newareabound}, we must show why these two horoball
arrangements are the critical configurations.

\noindent {\sl Proof of Proposition \ref{pro:newareabound}.} As above,
we may assume that, for all but finitely many values of $h$, $G \cap Q_h$
is a finite collection of immersed arcs.
The area of $G$ is at least
$$\int_0^\infty {{\rm Length}(G \cap Q_h) \over h} \, dh.$$
Since the aim of the proposition is to find a lower bound on the
area of $G$, we therefore will bound the length of
$G \cap Q_h$ from below by a function of $h$. Let us
first consider the case when $h > 1$. Then
$Q_h = \partial E_h$ is a torus. The surface $G \cap E_h$
forms a homology between $G \cap Q_h$ and
$\pm k [s]$ in $\partial M$. But any collection of curves in $Q_h$
homologous to $\pm k [s]$ must have length at least
$k \, {\rm Length}(s) /h$. 

Let us now focus on the case where $0 < h \leq 1$.
We define a function
$${\rm RelLength}(e_2,h) = \begin{cases}
1 & \text{if $h \geq 1$;} \\
\max \left \{ 0, 1 - {2 \over e_2} \sqrt{ 1 - h^2 } \right \} & \text{if $\min \{ 1/\sqrt 2, 1/e_2 \} \leq h \leq 1$;} \\
0 & \text{otherwise.}
\end{cases}$$
Geometrically, this is the ratio of the lengths of $\beta$
to $\alpha$ in Figure \ref{fig:rellen1} (if $e_2 \leq \sqrt 2$) and Figure \ref{fig:e2morethanroot2}
(if $e_2 > \sqrt 2$). Hence,
$$\int_0^\infty {e_ 2 \, {\rm RelLength}(e_2,h) \over h^2} dh$$
is the area of the shaded region in Figure \ref{fig:rellen1} or
Figure \ref{fig:e2morethanroot2}, which is $I(e_2) e_2$.

\noindent {\sl Claim 1.} The length of $G \cap Q_h$ is at least
$$k \, {\rm Length}(s)
{{\rm RelLength}(e_2,h) \over h}.$$

Thus, the claim asserts that, when $e_2 \leq \sqrt{2}$, the critical configuration is
shown in Figure \ref{fig:rellen1}, whereas when $e_2 > \sqrt{2}$, the critical configuration is
shown in Figure \ref{fig:e2morethanroot2}.

Let us assume the claim for a moment.
Then, the area of $G$ is at least
\begin{align*}
\int_0^\infty {{\rm Length}(G \cap Q_h) \over h} \, dh &\geq 
{k \, {\rm Length}(s) \over e_2} \int_0^\infty 
{e_ 2 \, {\rm RelLength}(e_2,h) \over h^2} dh \\
&= k \, {\rm Length}(s) I(e_2),\\
\end{align*}
thereby proving Proposition 4.3.

It is convenient to
give $\partial E_h$ the metric pulled back via the vertical
projection $\partial E_h \rightarrow \{ z=h \}$.
It then becomes a Euclidean torus.

The arcs $G \cap Q_h$ extend to a collection of closed curves
$G \cap \partial E_h$. The surface $G \cap E_h$
forms a homology between $G \cap \partial E_h$ and
$\pm k [s]$ in $\partial M$. Hence, the length of $G \cap \partial E_h$
is at least $k \, {\rm Length}(s) / h$. We wish to bound from below
the length of the parts of $G \cap \partial E_h$ that lie in $Q_h$.

We will shortly homotope $G \cap \partial E_h$ in
$\partial E_h$, creating new curves $C_h$. We will ensure that the length of
$C_h \cap Q_h$ is at most that of $G \cap Q_h$.
Thus, if we can show that $C_h \cap Q_h$ satisfies
the required lower bound on length, the claim will
be proved.  

\begin{figure}[htp]
\begin{center}
\includegraphics{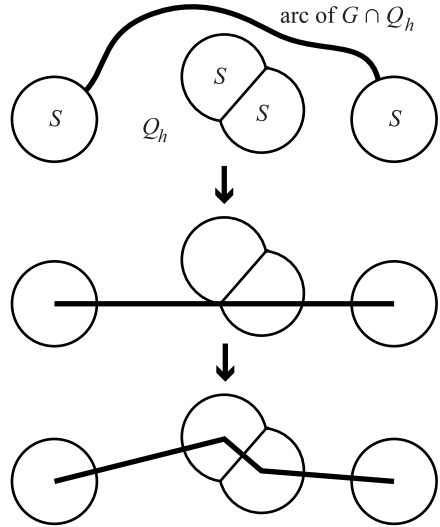}
\end{center}
\caption{Straightening the arcs $G \cap Q_h$}
\label{fig:arcstrgt}
\end{figure}

So, consider an arc of $G \cap Q_h$. We will now extend it to a curve which
sits in the Euclidean torus $\partial E_h$. The endpoints of the arc lie on totally geodesic faces of
$S$. Replace the arc by the geodesic arc that runs between
the points on these faces that are closest to $B_\infty$. This new arc may intersect
new faces of $S$ that the old arc did not.
If so, repeat this procedure. (See Figure \ref{fig:arcstrgt}.) The resulting curves
$C_h$ are a concatenation of Euclidean geodesic arcs.
The intersection of each such arc $\alpha$
with $S$ is either all of $\alpha$ or two sub-intervals of $\alpha$, each
of which contains an endpoint of $\alpha$. 

Since the curves $C_h$
are homologous to $\pm k [s]$, we obtain the inequality
$$\sum_\alpha {\rm Length}(\alpha) = {\rm Length}(C_h) \geq {k \, {\rm Length}(s) \over h}.$$

\begin{figure}[htp]
\begin{center}
\includegraphics{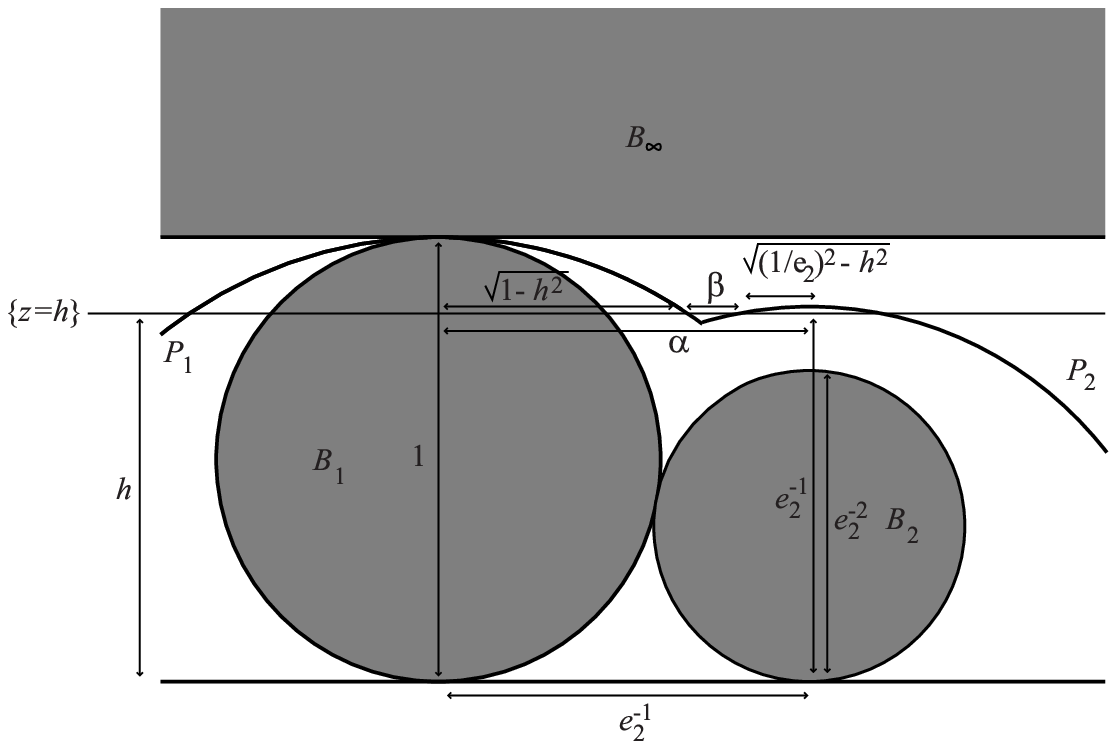}
\end{center}
\caption{}
\label{fig:rellen2}
\end{figure}

We now define a new function:
$${\rm RelLength}_2(e_2,h) = 
\max \left \{ 0, 1 - e_2\sqrt{ 1 - h^2 }  - \sqrt{ \max \{ 0, 1-e_2^2 h^2 \} } \right \}.$$
Geometrically, this is the ratio of the lengths of $\beta$
to $\alpha$ in Figure \ref{fig:rellen2}.

\noindent {\sl Claim 2.} When $e_2 \leq \sqrt{2}$,
$${{\rm Length}(\alpha \cap Q_h) \over {\rm Length}(\alpha)} \geq 
\min \{ {\rm RelLength}(e_2,h), {\rm RelLength}_2(e_2,h) \}.$$
When $e_2 > \sqrt{2}$,
$${{\rm Length}(\alpha \cap Q_h) \over {\rm Length}(\alpha)} \geq 
{\rm RelLength}(e_2,h).$$

Let $P_1$ and $P_2$ be the faces of $\tilde S$
containing the endpoints of a lift of $\alpha$.
These are equidistant planes between $B_\infty$
and, respectively, $B_1$ and $B_2$. 

Let us consider the case where $e_2 > \sqrt{2}$ first, as here the argument is much
simpler. Because $e_2 > \sqrt 2$, ${\rm RelLength}(e_2, h)$ is identically zero when $h \leq 1/e_2$.
Thus, we may assume that $h > 1/e_2$. The faces $P_1$ and $P_2$ contain
the endpoints of the lift of $\alpha$, and therefore they intersect the horosphere
$\{ z = 1/e_2 \}$. The
planes equidistant between $B_\infty$ and any horoball other than an ${\mathcal O}(1)$
horoball do not reach as high as $\{ z = 1/e_2 \}$. Thus,
$B_1$ and $B_2$ must be ${\mathcal O}(1)$ horoballs. The distance in $\partial B_\infty$
between their orthocentres is therefore at least $e_2$. Hence, it is clear that 
the ratio of the lengths of $\alpha \cap Q_h$ and $\alpha$ is at
least that of Figure \ref{fig:e2morethanroot2}, which proves Claim 2 when
$e_2 > \sqrt{2}$.

Thus, we now assume that $e_2 \leq \sqrt{2}$. If $B_1$ and $B_2$ are both
${\mathcal O}(1)$ horoballs, then, as above, the ratio of the lengths of
$\alpha \cap Q_h$ and $\alpha$ is at least that of Figure \ref{fig:rellen1},
which proves the claim in this case. Thus, we may assume that 
at least one of the horoballs
($B_2$, say) is not an ${\mathcal O}(1)$ horoball.
Let us suppose, for the sake of being definite,
that the distance between $B_1$ and $B_\infty$ is
no more than that between $B_2$ and $B_\infty$.
We will shortly perform some modifications to $B_1$ and $B_2$.
We will maintain each plane $P_i$ as the equidistant
plane between $B_i$ and $B_\infty$. Thus,
$P_i$ will be modified too. Each of these moves will
not increase ${\rm Length}(\alpha \cap Q_h)/{\rm Length}(\alpha)$. 
So, if we can show that the required
lower bound on ${\rm Length}(\alpha \cap Q_h)/{\rm Length}(\alpha)$
holds after these modifications, then it also held before.
In addition, each of these modifications will not decrease the radii
of $B_1$ and $B_2$. Hence, it will remain the case that the endpoints
of $\alpha$ do not lie in $Q_h$.

Translate $B_1$ and $B_2$ towards each
other until they become tangent. This reduces the lengths of $\alpha \cap Q_h$ and
$\alpha$ by the same amount, and so does not increase the ratio
${\rm Length}(\alpha \cap Q_h)/{\rm Length}(\alpha)$.
We next slide $B_1$ along $B_2$,
keeping them both just touching and moving the point
$B_1 \cap B_2$ closer to $B_\infty$. This has the effect
of enlarging $B_1$ and moving the orthocentre of $B_1$
away from the orthocentre of $B_2$. Thus, ${\rm Length}(\alpha)$
increases. Also, the union of $B_1 \cup B_2$
and everything below $B_1 \cup B_2$ has increased. Thus,
the subset of $P_1 \cup P_2$ consisting of points equidistant from
$B_\infty$ and $B_1 \cup B_2$,
has moved upwards. In particular, the old $\alpha \cap Q_h$
contains the new $\alpha \cap Q_h$, and therefore the length
of $\alpha \cap Q_h$ has not increased. Note that we are
using here the fact that the endpoints of $\alpha$ do not
lie in $Q_h$. Stop enlarging $B_1$ when it becomes tangent to $B_\infty$.

Next perform a similar slide, but with the roles of
$B_1$ and $B_2$ reversed, until the Euclidean diameter of $B_2$ is $e_2^{-2}$. 
Then, the horoballs are as shown in Figure \ref{fig:rellen2}, and so
$${{\rm Length}(\alpha \cap Q_h) \over {\rm Length}(\alpha)} \geq 
{\rm RelLength}_2(e_2,h).$$
Thus, we have proved Claim 2.

Claim 1 quickly follows from Claim 2 when $e_2 > \sqrt 2$, because
\begin{align*}
{\rm Length}(G \cap Q_h) & \geq \sum_\alpha {\rm Length}(\alpha \cap Q_h) \\
& \geq {\rm RelLength}(e_2,h) \sum_\alpha {\rm Length}(\alpha) \\
& \geq {\rm RelLength}(e_2,h)  k \ {\rm Length}(s) / h,
\end{align*}
as required. To prove Claim 1 when $e_2 \leq \sqrt 2$, we must compare
${\rm RelLength}(e_2,h)$ and ${\rm RelLength}_2(e_2,h)$.
Returning to the definitions of these functions, we see that we 
must compare
$$
{2 \over e_2}\sqrt{ 1-h^2} \quad \hbox{ and } \quad e_2 \sqrt{1-h^2} + \sqrt{\max \{ 0, 1 - e_2^2h^2 \}}.
$$
When $1 - e_2^2h^2 \leq 0$, the former is larger because $2/e_2 \geq e_2$. When
$1 - e_2^2h^2 > 0$,
\begin{align*}
&\ \quad {2 \over e_2}\sqrt{ 1-h^2} \geq e_2 \sqrt{1-h^2} + \sqrt{1 - e_2^2h^2} \\
&\Leftrightarrow ( {2 \over e_2} - e_2 ) \sqrt{1-h^2} \geq \sqrt{ 1 - e_2^2h^2} \\
&\Leftrightarrow \left ( {4 \over e_2^2} - 4 + e_2^2 \right ) (1-h^2) \geq
1 - e_2^2 h^2 \\
&\Leftrightarrow 
{4 \over e_2^2} - 5 + e_2^2 - {4 h^2 \over e_2^2} + 4h^2 \geq 0 \\
&\Leftrightarrow 4h^2 \left( 1- {1 \over e_2^2} \right ) \geq
\left ( {4 \over e_2^2} - 1 \right ) (e_2^2 - 1) \\
&\Leftrightarrow h^2 \geq \left ( 1 - {e_2^2 \over 4} \right ).
\end{align*}
Thus, when $h^2 \geq  1 - (e_2^2 /4)$, 
${\rm RelLength}_2(e_2,h) \geq {\rm RelLength}(e_2,h)$.
Moreover, we have equality when $h^2 =  1 - (e_2^2 /4)$.
But when $h^2 = 1 - (e_2^2 /4)$, we see from Figure \ref{fig:rellen1}
that $\alpha \cap Q_h$ is a single point and so
${\rm RelLength}(e_2,h)$ is zero. Hence, ${\rm RelLength}_2(e_2,h)$
is also zero. Since ${\rm RelLength}$
and ${\rm RelLength}_2$ are non-decreasing functions of
$h$, we deduce that they are both zero when 
$h^2 < 1 - (e_2^2 /4)$. Thus, ${\rm RelLength}(e_2,h)$ is,
for all values of $h$, the minimum of the two quantities.

This proves Claim 1 and hence the proposition. \qed\

\section{Area control}

As discussed in Section 2,
if we can get good lower bounds for the area of the cusp torus of $M = {\Bbb H}^3/\Gamma$, then we will be able to 
fruitfully control the number of possible exceptional slopes.  
We now begin developing this area control by analyzing the maximal cusp diagram for $M$.

Consider our standard, normalized lift of the cusp neighborhood to upper-half-space.
The view from infinity consists of a collection of overlapping disks.  
Specifically, we vertically project all horoballs other than $B_\infty$ to the plane $z = 1$, 
thereby producing a collection of disks, the largest of which have radius one-half.  
Let $\sigma$ be a parabolic transformation in $\Gamma$ that preserves $B_\infty$
and has minimal translation length.  We may rotate the picture
about the $z-$axis so that $\sigma$ takes $(0,0)$ to $(m,0)$ where $m > 0$. 
Note that, automatically, $m \ge 1$.  We take $\tau$ to be another
element of $\Gamma$, so that $\sigma$ and $\tau$ together generate
the stabiliser of $B_\infty$. We can assume that $\tau$ 
takes $(0,0)$ to $(tm, h)$ where $1/2 < t \leq 1/2$ and $h > 0$. 
The {\sl maximal cusp diagram} ${\mathcal M}$ for $M$ consists of 
this euclidean plane with all the vertically projected disks, 
and arrows representing $\sigma$ and $\tau$. Quotienting
the maximal cusp diagram by the stabiliser of $B_\infty$ produces
the cusp torus.

Each orthoclass ${\mathcal O}(i)$ is an equivalence class of horoballs.
The image of these horoballs in ${\mathcal M}$ is a collection of disks.
They come in two orbits under the action of 
the stabiliser of $B_\infty$. (The fact that there
are two orbits was first observed by Adams \cite{Ada1}.
One of the associated horoballs is known as an {\sl Adams} horoball.)
Thus, their image in the cusp torus is two discs, which we denote $D_i$ and $D_i'$.
We also refer to the disks in
${\mathcal M}$ as ${\mathcal O}(i)$ disks.

The disks $D_1$ and $D_1'$ each have radius one-half and so they contribute 
$2\pi (1/2)^2 = \pi/2$ to the area of the cusp torus.

We now want to use the next largest disks $D_2$ and $D_2'$ in the cusp torus to get more area. 
We observe that if $o(2)$ is small (that is, close to 0) then 
the $D_2$ and $D_2'$ disks are large (radius close to one-half), while if $o(2)$ is large, 
then the $D_2$  and $D_2'$ disks are small.  

By Lemma \ref{lem:orthocentreseparation}, if an ${\mathcal O}(i)$ disk 
and an ${\mathcal O}(j)$ disk have associated horoballs separated by $o(k)$ then the distance between their centers is 
$e_k/(e_i e_j)$. By Lemma \ref{lem:nonnn}, any two ${\mathcal O}(1)$ disks in ${\mathcal M}$ must have centers 
separated by at least $e_2$.  Hence, we can extend the $D(1)$ and $D(1)'$ disks to have radius $e_2/2$ and 
still have their interiors be disjoint. So, $2 \pi (e_2/2)^2 = \pi e_2^2/2$ is a lower bound for the area 
of the cusp torus. Let $D_1^+$ denote the union of these two enlarged disks. 

Further, we can add on the area provided by the $D_2$ and $D_2'$ disks.  These 
disks have radius $1/(2 e_2^2)$. But we need to take into account the fact that the $D_1^+$ 
and $D_2/D_2'$ disks might overlap. The largest overlap occurs when the associated 
horoballs are abutting; in which case their centers are a distance $e_1/(e_1 e_2) = 1/e_2$ apart. 
 Another problem is that there may be more than one overlap.  For example, a $D_2/D_2'$ disk might be overlapped 
by both $D_1^+$ disks.  

In order to control the number of overlaps, we use Mom Technology (see Section 3). 
Loosely, the theme is that too much overlap leads to Mom structures. 
Note that in the case of no overlaps, we have the nice situation 
that when $e_2$ is ``large" we get a big area contribution from $2 \pi(e_2/2)^2$ alone, and in the 
``small" $e_2$ situation the $D_2$ and $D_2'$ disks are big and so we still get a ``big" area.

Assuming that $M$ does not contain a geometric Mom-2 involving ${\mathcal O}(1)$ and ${\mathcal O}(2)$, 
it turns out that we can quickly improve our previous estimate, as follows. Simply 
expand the $D_1^+$ disks to radius $e_3/2$ (not simply radius $e_2/2$).  This will result 
in overlap if the centers of two ${\mathcal O}(1)$ disks are within $e_2$ of each other. 
But the overlap can be controlled because we are in the no-Mom-2 situation.  Further,
we put disks $D_2^+$ of radius $e_3/e_2 - e_3/2$ centered at the $D_2$ and $D_2'$ disks (this radius is chosen 
so as to give fruitful area, but to avoid overlap not handled by the no-Mom-2 condition).
The no-Mom-2 condition means that we have to consider only two types of overlap cases.  

First, the $(1,1,2)$ case where there are three horoballs $B_1,\ B_2,\ B_3$ 
which form a $(1,1,2)$-triple.  In ${\mathcal M}$, this manifests itself as 3 overlaps, one for each of 
the 3 horoballs being sent to $B_\infty$ by an element of $\Gamma$.  Thus, because $B_1$ is $o(1)-$separated 
from $B_2$ and $B_3$, then mapping $B_1$ to $B_\infty$ results 
in two ${\mathcal O}(1)$ disks whose centers are a distance $e_2$ apart.  This results in an overlap because the 
expanded disks have radius $e_3/2$.  Sending $B_2$ to $B_\infty$ results in an ${\mathcal O}(1)$ disk and an 
${\mathcal O}(2)$ disk whose centers are separated by $e_1/(e_1 e_2) = 1/e_2$ and again there is overlap for 
the expanded disks at these centers.  Sending $B_3$ to $B_\infty$ produces the same overlap picture.

However, because we are assuming that $M$ does not have a geometric Mom-2 involving
${\mathcal O}(1)$ and ${\mathcal O}(2)$,
there are no other overlaps between the $D_1^+$ and $D_2^+$ discs.
Thus, in the case where there is a $(1,1,2)$-triple, we get the following area lower bound for the cusp torus:
\begin{align*}
& 2 \pi (e_3/2)^2 + 2 \pi \left ({e_3 \over e_2} - {e_3 \over 2} \right)^2 - 
{\rm overlap}(e_3/2, e_3/2, e_2)\\
& \qquad - 2 \ {\rm overlap}\left ({e_3 \over 2}, {e_3 \over e_2} - {e_3 \over 2}, {1 \over e_2} \right)
\end{align*}
 where ${\rm overlap}(a, b, c)$ is the area of the overlap of two disks of radius 
$a$ and $b$ whose centers are separated by $c$.

The second case is when there is a $(1,2,2)$-triple and here the lower bound for area is:
\begin{align*}
& 2 \pi (e_3/2)^2 + 2 \pi \left ({e_3 \over e_2} - {e_3 \over 2} \right)^2 - 
2 \ {\rm overlap}\left({e_3 \over 2}, {e_3 \over e_2} - {e_3 \over 2}, 1 \right) \\
& \qquad - {\rm overlap}\left({e_3 \over e_2} - {e_3 \over 2}, {e_3 \over e_2} - {e_3 \over 2}, 
{1 \over e_2^2} \right)
\end{align*}
By analyzing the overlap function, it can be seen that to find a valid lower bound only the 
$(1,1,2)$ case is needed.  In particular, we exploit the following result.

\begin{lemma}
If $a_1 \ge  a_2,\ b_1 \ge b_2$ and $a_1 + b_1 - c_1 \ge a_2 + b_2 - c_2$ 
then $${\rm overlap}(a_1, b_1, c_1) \ge {\rm overlap}(a_2, b_2, c_2).$$
\label{lem:overlap}
\end{lemma}

\noindent{\sl Proof.}
$a_1 + b_1 - c_1$ is the {\sl linear dimension of overlap}, 
that is, it's the length of the intersection of the physical overlap 
with the line connecting the centers.  So, the first pair of circles 
(the circles of radius $a_1$ and $b_1$) have larger radii and a greater linear dimension 
of overlap than the second pair of circles.  Hence, the area of the overlap is larger for 
the first pair of circles than for the second pair of circles.  The point here is that not 
only is the overlap wider in the first case (linear dimension of overlap is greater) but also, 
in the perpendicular direction (to the linear dimension of overlap) the boundaries of the overlap 
are more vertical because the associated circles have larger radii. \qed\

In comparing the $(1,1,2)$ and $(1,2,2)$ cases, there are 3 overlap comparisons.
First compare $112$ versus $122$.  That is, compare the two $D_1^+$
disks whose centers are $e_2$ apart with a $D_1^+$ and $D_2^+$ disk
whose centers are $e_2/(e_1 e_2) = 1$ apart.  An application of Lemma \ref{lem:overlap}
shows the overlap 
contribution is greater for $112$.  Next compare $121$ and $221$, and then compare $211$ and $221$. 
In both cases, the overlap is greater for the first element in the comparison.  Hence, the total overlap 
punishment is larger in the $(1,1,2)$ case.

Thus, we have the following result.

\begin{theorem}
Let $M$ be a compact orientable 3-manifold, with boundary
a torus, and with interior admitting a complete finite-volume hyperbolic structure. Suppose that $M$ does not contain
a geometric Mom-2 involving ${\mathcal O}(1)$ and ${\mathcal O}(2)$. 
Then, the area of the cusp torus is at least
\begin{align*}
& 2 \pi (e_3/2)^2 + 2 \pi \left ({e_3 \over e_2} - {e_3 \over 2} \right )^2 - 
{\rm overlap}(e_3/2, e_3/2, e_2)\\
& \qquad - 2 \ {\rm overlap}\left ({e_3 \over 2}, {e_3 \over e_2} - {e_3 \over 2}, {1 \over e_2} \right).
\end{align*}
\label{thm:areanomom2}
\end{theorem}

Just using the $D_1^+$ and $D_2^+$ disks does not provide enough area for our purposes. 
 So we analyze the ${\mathcal O}(3)$ disks, use $e_4$, and exploit Mom-3 technology.  That is, we expand 
the $D_1^+$ disks to radius $e_4/2$, we expand the  $D_2^+$ disks to radius 
$e_4/e_2 - e_4/2$, and we expand the $D_3$ and $D_3'$ disks to radius 
$e_4/e_3 - e_4/2$, giving disks $D_3^+$.  To control overlap, we assume that $M$ does not have a 
torus-friendly geometric Mom-3 involving ${\mathcal O}(1)$, ${\mathcal O}(2)$ and ${\mathcal O}(3)$. 

There are a variety of cases where overlaps do not yield a torus-friendly geometric Mom-3 (see \cite{GMM3} for the list) 
and we need to do overlap analysis for these. We must consider the overlap
generated by an $(a,b,c)$ triple, where $a,b,c \in \{ 1,2,3 \}$.
As in the $(1,1,2)$ case, there are (at most) three overlaps, with area
\begin{align*}
& {\rm overlap}(R(a), R(b), e_c/(e_a e_b)),\\
& {\rm overlap}(R(b), R(c), e_a/(e_b e_c)),\\
& {\rm overlap}(R(c), R(a), e_b/(e_c e_a)),
\end{align*}
where $R(a)$, $R(b)$ and $R(c)$ are the radii of the discs $D_a^+$, $D_b^+$ and $D_c^+$.
We now use the following:

\begin{lemma}
For integers $a,b,c,d,e,f \in \{ 1,2,3\}$ such that $a \geq d$, $b \geq e$ and $c \geq f$,
$${\rm overlap}(R(a),R(b),e_c/(e_a e_b)) \leq {\rm overlap}(R(d),R(e),e_f/(e_d e_e)).$$
\label{lem:overlap2}
\end{lemma}

\noindent {\sl Proof.} Note first that, in order to prove the lemma, we may
assume that two of the inequalities $a \geq d$, $b \geq e$ and $c \geq f$ are
actually equalities. Note also that $R(x) = (e_4 / e_x) - (e_4/2)$. Hence,
$R(a) \leq R(d)$ and $R(b) \leq R(e)$. In order to apply Lemma
\ref{lem:overlap}, we therefore need to know that 
$$R(a) + R(b) - e_c/(e_a e_b) \leq R(d) + R(e) - e_f / (e_d e_e).$$
This is clear if $a = d$ and $b = e$. If $b = e$ and $c = f$, the
inequality becomes 
$$e_4 / e_a  - e_4 / 2 - e_c / (e_a e_b) \leq e_4 / e_d - e_4/2 - e_c / (e_d e_b),$$
which is equivalent to 
$${e_4 - (e_c / e_b) \over e_a} \leq {e_4 - (e_c / e_b) \over e_d},$$
and this holds because $e_a \geq e_d$ and $e_4 \geq e_c / e_b$.
\qed

Thus, we see that an $(a,b,c)$-triple will produce less overlap area than a $(d,e,f)$-triple if 
$a \ge d,\ b \ge e,\ c \ge f$.  Exploiting this observation, there are two possible
cases that result in the maximum total overlap area: $(1,1,2), (1,1,3)$ and
$(1,2,3),(1,2,3),(1,1,2)$. This gives the following result.

\begin{theorem}
Let $M$ be a compact orientable 3-manifold, with boundary
a torus and with interior admitting a complete finite-volume hyperbolic structure. 
Suppose that $M$ does not contain
a torus-friendly geometric Mom-3 involving ${\mathcal O}(1)$, ${\mathcal O}(2)$ and ${\mathcal O}(3)$,
or a geometric Mom-2 involving a subset of ${\mathcal O}(1)$, ${\mathcal O}(2)$ and ${\mathcal O}(3)$.
Then, the area of the cusp torus is at least the minimum of
\begin{align*}
& 2 \pi (e_4/2)^2 + 2 \pi \left({e_4 \over e_2} - {e_4 \over 2} \right)^2 
+ 2 \pi \left( {e_4 \over e_3} - {e_4 \over 2} \right)^2 \\
& \qquad - {\rm overlap}(e_4/2,e_4/2,e_2) - {\rm overlap}(e_4/2,e_4/2,e_3) \\
& \qquad - 2 \ {\rm overlap}(e_4/2, e_4/ e_2 - e_4/2, 1/e_2) \\
& \qquad - 2 \ {\rm overlap}(e_4/2, e_4/e_3 - e_4/2, 1/e_3)
\end{align*}
and 
\begin{align*}
& 2 \pi (e_4/2)^2 + 2 \pi \left({e_4 \over e_2} - {e_4 \over 2} \right)^2 
+ 2 \pi \left( {e_4 \over e_3} - {e_4 \over 2} \right)^2 \\
& \qquad - 2 \ {\rm overlap}(e_4/2, e_4/ e_2 - e_4/2, e_3/e_2) \\
& \qquad - 2 \ {\rm overlap}(e_4/e_2 - e_4/2, e_4/ e_3 - e_4/2, 1/(e_2e_3)) \\
& \qquad - 2 \ {\rm overlap}(e_4/e_3 - e_4/2, e_4/2, e_2/e_3) \\
& \qquad - 2 \ {\rm overlap}(e_4/2, e_4/ e_2 - e_4/2, 1/e_2) \\
& \qquad - {\rm overlap}(e_4/2,e_4/2,e_2).
\end{align*}
\label{thm:areanomom3}
\end{theorem}

In fact, under certain circumstances, it is possible to improve this yet further.
We place disks of radius $1/(e_4e_2) - e_4/e_2 + e_4/2$ at the centers of
$D_4$ and $D_4'$. It is shown in \cite{GMM3} that these disks
do not overlap with the $D_2^+$ and $D_3^+$ disks. Moreover,
if $2e_4^2 + 2e_4 - e_2(e_4^2 + e_4 + 1) \geq 0$, then these
disks are embedded and disjoint.
They might overlap with the $D_1^+$ disks. However, if $e_2 + 1 \geq e_4^2$,
then this triggers a $(1,1,4)$-triple. Thus, assuming in addition
that $M$ does not have a geometric Mom-2 involving ${\mathcal O}(1)$
and ${\mathcal O}(4)$, there can only be one such overlap.
Hence, we obtain the following. (See Lemma 5.6 in \cite{GMM3}.)

\begin{theorem}
Let $M$ be a compact orientable 3-manifold, with boundary
a torus, and with interior admitting a complete finite-volume hyperbolic
structure. Suppose that $M$ contains neither a torus-friendly geometric Mom-3
involving ${\mathcal O}(1)$, ${\mathcal O}(2)$ and ${\mathcal O}(3)$ nor
a geometric Mom-2 involving a subset of ${\mathcal O}(1)$, ${\mathcal O}(2)$, ${\mathcal O}(3)$ and ${\mathcal O}(4)$.
Suppose also that $2e_4^2 + 2e_4 - e_2(e_4^2 + e_4 + 1) \geq 0$
and that $e_2 + 1 \geq e_4^2$. 
Then, the area estimate in Theorem \ref{thm:areanomom3}
can be increased by
$$
2 \pi \left ( {1 \over e_2 e_4 }- {e_4 \over e_2} + {e_4 \over 2} \right)^2
- 2 \ {\rm overlap}\left ({e_4 \over 2}, {1 \over e_4 e_2} - {e_4 \over e_2} + {e_4 \over 2}, {1 \over e_4} \right).
$$
\label{thm:bonusball}
\end{theorem}

Note that when $e_4 < (1 + \sqrt 5)/2$, the condition $2e_4^2 + 2e_4 - e_2(e_4^2 + e_4 + 1) \geq 0$
is satisfied. This is because $e_2 \leq e_4$ and so 
$$2e_4^2 + 2e_4 - e_2(e_4^2 + e_4 + 1) \geq -e_4^3 + e_4^2 + e_4 > 0.$$

\section{Tools for the parameter space analysis}

In this section, we describe our tools for excluding regions of the parameter space.
Recall from Section 1 that we are using 6 parameters: $e_2$, $e_3$,
$e_4$, $m$, $t$ and $h$ and here, we work as though these parameters
are fixed and given. However, in practice, they are only specified to lie
within small intervals. The parameters $m$, $t$ and $h$ specify a lattice,
which consists of the orthocentre of an ${\mathcal O}(1)$ horoball and its images
under the covering transformations that preserve $B_\infty$.
One lattice point is at the origin $(0,0)$. A closest lattice point to $(0,0)$ is
at $A = (m,0)$. A closest lattice point with non-zero second co-ordinate
lies at $B = (mt,h)$, where $-{1 \over 2} \leq t \leq {1 \over 2}$
and $h > 0$. In fact, reflecting the picture if necessary, we may assume
that $0 \leq t \leq {1 \over 2}$.

We warm up by analyzing a few representative parameter points.

First, $e_2 = 2.0,\ m = 2.0, \ h = 3.2$.  The fact that $e_2 = 2.0$ tells us that the centers of 
the full-sized disks in the maximal cusp diagram ${\mathcal M}$ are separated by at least distance $2.0$.  
As above,  we can embed 2 disks of radius 1 in ${\mathcal M}$.  These contribute area $2 \pi$ to the area 
of the cusp torus.  By disk-packing, we can improve this lower bound on area to $(2 \pi) (2 \sqrt{3}/\pi) = 4 \sqrt{3}.$
Hence, no parameter point with these values for $e_2, m, h$ can be realized because the $e_2$ value implies 
that the area of the cusp torus is at least $4 \sqrt{3}$ but the values of $m$ and $h$ imply that the area of 
the cusp torus is $6.4$, a contradiction. 

Second, $e_2 = 2.0,\ m = 2.0, \ h = 4.0$.  There is no immediate area contradiction here, so we analyze 
slopes in ${\mathcal M}$.  By Theorem \ref{thm:improve6}, we know the only possible exceptional slopes $(p,q)$ must have 
associated lattice point $(pm + qtm, qh)$ that is within $L = 6 \pi / (\pi + 6 - 3 \sqrt 3)$ of the origin.  
That is, we need to have  $(p + qt)^2m^2 + q^2h^2 \le L^2 \le 23$, which becomes 
$(p + qt)^2 \le (23 - q^2 h^2)/(m^2).$  Because $h = 4$ in our example, we can see that there can be 
no exceptional slopes when $q \ge 2.$  When $q = 1$ our equation becomes
$(p + qt)^2 \le (23 - 16)/(m^2) = 1.75$.  So, when $q = 1$, we have 
$|p + qt| \le 1.33$ and this can occur for at most 3 points.
Together with $(1,0)$, the unique $q = 0$ slope, we end up
with at most $4$ exceptional slopes for these parameter values.
This holds regardless of the values of $e_3,\ e_4,\ t$, and further, it holds
when $h > 4$ because then the inequalities work at least as well.

Third, 
$e_2 = 1.26,\ 
e_3 = 1.38,\ 
e_4 = 1.38,\ 
m = 2.19463$.
Ignoring $t$ and $h$ for the time being, we can compute a lower bound
for the area of the cusp torus; using Theorems \ref{thm:areanomom3} and \ref{thm:bonusball}, 
we get that the area is at least  $4.13103.$
As above, the formula for possible exceptional slopes $(p,q)$ is
$(p + qt)^2m^2 + q^2h^2 \le L^2$ where $L$ is determined as in Theorem \ref{thm:improve6}.
Here we will take $h = 4.13103/m$ (if $h$ is less than this number 
then there is an immediate contradiction).

In particular, when $q = 1$ we have 
$|p + tq| < 2.50441$.   When $t = 0$ this means there are most 5 exceptional
slopes with $q = 1$, but when $t = 0.5$ we see that there could be 6 exceptional slopes.
For $q = 2$ we get $|p + tq| < 2.01622$.  When $t = 0$ this means there are most 2 exceptional
slopes with $q = 2$ (note that the $(p,q)$ must be relatively prime pairs of integers),
 but when $t = 0.5$ we see that there could be 3 exceptional slopes.
For $q = 3$ we get $|p + tq| < 0.62201$.  When $t = 0$ this means there are no exceptional
slopes with $q = 3$, but when $t = 0.5$ we see that there could be 2 exceptional slopes.
Adding in the exceptional slope $(1,0)$ yields a maximum of 8 exceptional slopes when $t = 0.$
Further, the analysis holds for larger values of $h$, so more parameter points are 
eliminated for free.

But when $t = 0.5$ we have a maximum of 12 exceptional slopes and the associated
parameter point is not yet eliminated.  So, we do the following simple check:
Consider the Adams horoball, and determine how close its orthocenter $C$ is to the vertices 
of the triangle determined by the origin $O$, the point $A = (m,0)$ and the point $B = (tm,h).$
It turns out that there is no point in the triangle which is simultaneously further
than $1.26101$ from the vertices (1.26101 is the circumradius of the triangle $OAB$).
If this number is less than $e_2$ then we have an immediate contradiction to the fact
that $e_2$ is the second shortest Euclideanized ortholength.  Unfortunately,
$e_2$ is just slightly less than the circumradius, and there is no immediate contradiction.
So, we have to do a slightly more complicated check.  Because we are assuming that
there is no geometric Mom-2 involving ${\mathcal O}(1)$ and ${\mathcal O}(2)$,
we see that in actuality, the Adams horoball orthocenter $C$ can be $e_2$ away
from at most one vertex, and that it is at least $e_3$ from the other two vertices. 
By a straightforward calculation, for the parameter point at hand, this can't happen 
and we have eliminated this parameter point.

In analyzing parameter points we have trivial tools like the one used in the first example above, and 
we have 3 non-trivial  tools:  First, compute a lower bound on area and hence an upper bound on the
number of exceptional slopes $(p,q)$.
Second, use circumradius.  Third, bring in the no-Mom-2 assumption to further constrain the possible 
position of the center of the Adams horoball.

More specifically, Tool 1 is as follows.

{\sl Tool 1}: For each $e_2$, $e_3$ and $e_4$, compute a lower bound on the area
of the cusp torus, using the formulae from Section 5. Then, for each torus
with at least this area, use the improved version of the 6-theorem to
bound the set of exceptional slopes.

We now describe Tools 2 and 3 in more detail. Tool 2 is used to show that
certain regions of the parameter space lead to a contradiction. Tool 3 is used
to show that, for certain other regions, $M$ contains a geometric Mom-2
involving ${\mathcal O}(1)$ and ${\mathcal O}(2)$.

Consider the orthocentre $C$ of the Adams horoball.
Its distance from each lattice point is at least $e_2$. In fact, if
it closer than $e_3$ to a lattice point, then this implies that 
there is a $(1,1,2)$-triple of horoballs. So, if it is closer than $e_3$ to 
two lattice points, there are two inequivalent $(1,1,2)$-triples and so $M$
contains a geometric Mom-2 involving ${\mathcal O}(1)$ and ${\mathcal O}(2)$.

{\sl Tool 2:} 
This consists of one test. Check whether there is some point that has distance
at least $e_2$ from all lattice points. If not, this arrangement cannot occur and 
such parameter points are eliminated. If so, then this arrangement is provisionally
permitted. Note that, here, no Mom technology is
required.

{\sl Tool 3:}
This consists of two tests:
\begin{enumerate}
\item Check whether there is some point that has distance at least $e_3$ from all lattice
points. If so, then this arrangement is permitted. If not, then
pass to the second test.
\item The orthocentre $C$ is therefore distance less than $e_3$ from
a lattice point, which we may take to be $(0,0)$. So, it lies
within an annulus centred at $(0,0)$, with inner radius $e_2$
and outer radius $e_3$. The second test verifies whether there
is a point in this annulus that is distance at least $e_3$ from
every other lattice point. If so, then this arrangement is permitted.
If not, then $M$ must contain a geometric Mom-2 involving ${\mathcal O}(1)$ and ${\mathcal O}(2)$.
\end{enumerate}

In order to implement  Tool 2 and Test 1 of Tool 3, we compute the
circumradius of the triangle $T$ with corners
$O = (0,0)$, $A = (m,0)$, $B = (mt,h)$, which is
$${|OA| \ |AB| \ |BO| \over 4 {\rm Area}(T)} = {m \sqrt { m^2 (1-t)^2 + h^2}
\sqrt { m^2 t^2 + h^2} \over 2 mh}.$$
Clearly, if this is definitely less than $e_2$, then every point in $T$
has distance less than $e_2$ from one of $O$, $A$ and $B$. Hence,
every point in the plane is less than $e_2$ from some lattice point, and
hence this parameter point can be discarded. This is Tool 2. If the
circumradius is at least $e_3$, then Test 1 of Tool 3 passes. Otherwise,
we pass to Test 2, which we now describe.

\begin{figure}[htp]
\begin{center}
\includegraphics{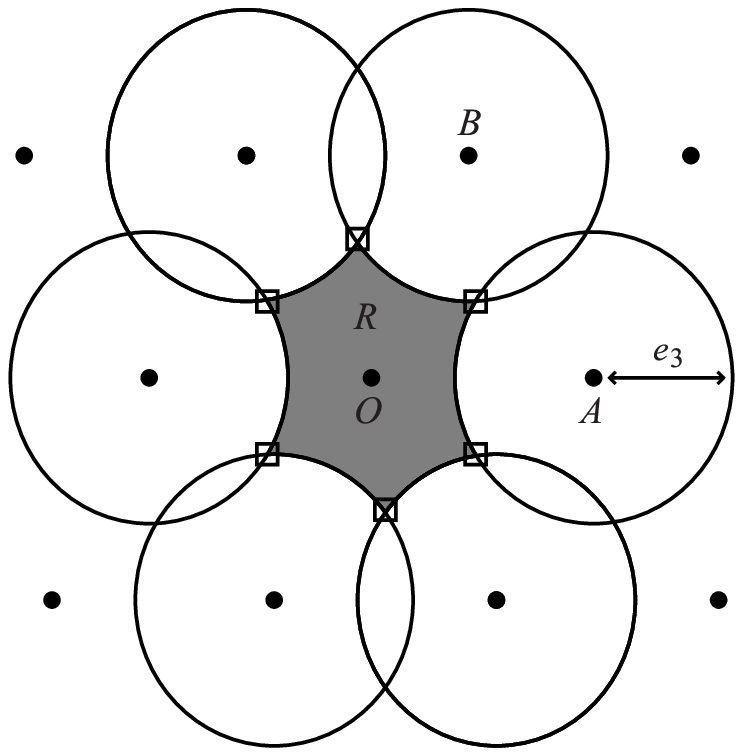}
\end{center}
\caption{}
\label{fig:hexacirc}
\end{figure}

Assuming that Test 1 of Tool 3 has failed, the circumradius of $T$ is less than
$e_3$. Thus, the discs of radius $e_3$ centred at the corners of $T$
cover $T$. Hence, the sides of $T$ each have length at most $2 e_3$.
(Note that the midpoint of each side is closer to the endpoints of
the side than the remaining vertex, because the triangle $T$ is not obtuse.)
So, the six primitive lattice points that are closest to $O$
form a hexagon with side lengths less than $2e_3$. Therefore,
the discs of radius $e_3$ centred at these points enclose
a region $R$, as shown in Figure \ref{fig:hexacirc}. We now define 6 distinguished points.
Consider one of the six copies of $T$ with $O$ as a vertex. Place two circles
of radius $e_3$ at its two vertices which are not $O$. Then
we consider the point of intersection between these two circles
that is closest to $O$. This is one of the 6 distinguished points.
In Figure \ref{fig:hexacirc}, these 6 distinguished points are marked with small squares.

\begin{lemma}
Suppose that Test 1 of Tool 3 has failed. If all 6 of these distinguished points
have distance less than $e_2$ from $O$, then $M$ contains a geometric Mom-2
involving ${\mathcal O}(1)$ and ${\mathcal O}(2)$.
\label{lem:hex}
\end{lemma}

\noindent {\sl Proof.} Consider one of the 6 primitive lattice points closest
to $O$. The circle of radius $e_3$ about this vertex contains two distinguished
points. We join these two points by the arc of the circle
that is shorter. The union of these arcs forms a closed curve $F$
shaped like a hexagon. The boundary of $R$ forms a subset of $F$.
Typically, the boundary of $R$ will be all of $F$, but there are
situations where it is not. For example,
if $|OA| < e_3$, then the discs centred at $A$ and its
reflection in the origin overlap. 

We claim that a point of $F$ with maximum distance from $O$ is at
one of the 6 distinguished points. This follows from the fact that
$F$ is a union of arcs of circles. For any such arc, a point of
maximal distance from $O$ is at an endpoint of the arc. This is
because the unique point of the circle at maximal distance from $O$
is not contained in the arc.

Thus, if all $6$ distinguished points have distance less than $e_2$ from $O$,
then this is true for all points of $F$, and hence all of $R$,
since $\partial R$ is a subset of $F$. 

Now, if there is no geometric Mom-2 involving ${\mathcal O}(1)$ and ${\mathcal O}(2)$, 
then $C$ is distance at least $e_3$
from every lattice point other than $O$. It therefore misses the interiors of the
discs, and so lies in $R$. But we have seen that every point of
$R$ has distance less than $e_2$ from $O$, and $C$ is not
permitted to lie this close to $O$. Thus, we deduce in fact that
$M$ has a geometric Mom-2 involving ${\mathcal O}(1)$ and ${\mathcal O}(2)$. \qed\

Thus, Test 2 of Tool 3 computes the distances of
these six distinguished points from $O$, and determines whether they are all less
than $e_2$. If so, the configuration
is not permitted. Otherwise, it is. In fact, only three
of these points need to be tested, since $R$ is symmetric under reflection
in the origin and so opposite points have the same distance
from $O$.

\section{The parameter space analysis}

Here we describe our parameter space analysis in detail.
Our first job is to reduce to a compact subset of the parameter
space.   

We begin by controlling $h, m, t.$

One crude reduction is to note that when $h > 4$ 
we can quickly show that no more than 10 exceptional slopes
are possible.  That is, when $h > 4$,
there are no exceptional slopes possible for $q \ge 2$ regardless
of the critical slope length $L(e_2)$ from Theorem \ref{thm:improve6}, and when $q = 1$ there can be no more that 
9 exceptional slopes.  Specifically, when $q = 1$, the  exceptional
slope analysis comes down to $(p + tq)^2 \le (L(e_2)^2 - h^2)/(m^2) < (36 - 16)/1 = 20$,
and we see that $|p + tq|$ must be less than $4.5$.  Hence, independent of the value
of $t$ there can be no more than 9 exceptional slopes when $q = 1$. So, 
there can be no more than 10 exceptional slopes and the intersection number
of any two such slopes is at most $8$.

In addition, $h$ must be at least $1.7$. For, the fact that $m$ is the
minimum translation length implies that $h \geq m \sqrt 3/2$. Cao and Meyerhoff \cite{CM1}
give that the cusp area is at least 3.35. Thus,
$3.35 \leq hm \leq 2 h^2 / \sqrt 3$, and hence $h > 1.7$.

Note that this implies that if $(p,q)$ is an exceptional slope, then
$|q| \leq 3$. This is because the length of the slope is at least
$|q|h$, which is more than $6$ when $|q| \geq 4$.

The minimum translation length $m$ must satisfy $m \ge e_2$, by Lemmas \ref{lem:nonnn}
and \ref{lem:orthocentreseparation}.
Further, $m > 2.5$ can be shown to result in at most 8 exceptional slopes. 
That is, $m > 2.5$ implies that $h > 2$ and hence there are no $|q| \ge 3$ 
exceptional slopes.  For $q = 2$ we get
$(p + tq)^2 \le (L(e_2)^2 - 4 h^2)/(m^2) < (36 - 16)/(6.25) < 4,$
hence $|p + tq| < 2$ and there are at most 2 exceptional slopes with $q = 2.$
Finally, when $q = 1$ we get 
$(p + tq)^2 \le (L(e_2)^2 -  h^2)/(m^2) < (36 - 4)/(6.25) < 6,$
hence $|p + tq| < 2.5$ and there are at most 5 exceptional slopes with $q = 1.$
It is also easy to check that in this case, the intersection number of any
two exceptional slopes is at most $5$.
So, we can restrict to $e_2 \le m \le 2.5.$

By symmetry, we can assume $0.0 \leq t \leq 0.5$.

Now we control $e_2$.

When $e_2$ is large (for example $e_2 = 2$) we can use the 
elementary area argument from Section 5 using disks with radius $e_2/2$ to get good area control.
In fact, we can have the computer analyze the parameter space with $m,\ h,\ t$ 
as restricted above and with $e_2 \leq 2.0$ (the $e_3$ and $e_4$ parameters
follow, for free) and show that at each point, the maximum number of 
exceptional slopes is at most 10 and the intersection number between
any two exceptional slopes is at most 8. This is easily done by using interval 
arithmetic to break up the parameter space into small sub-boxes
and then having the computer analyze each small sub-box.  
Tool 1 alone establishes the theorems as $e_2$ descends from
$2.0$ to $1.5$. Of course, for $e_2$ larger than 2 we get more area by this approach
and the tool 1 argument always works.
At a little less than $e_2 = 1.5$, the tool 1 analysis breaks down 
and we can start applying tool 2, 
the circumradius tool.  Utilizing this tool
enables us to get $e_2$ down to $1.4$.  This program is called {\tt slopes1}
and is available from the authors \cite{LM1}.
Thus, we can restrict to $1 \le e_2 \le 1.4$.

We now work on controlling $e_3$.  

We work under the assumption that $1 \le e_2 \le 1.4$.
The area control provided by our $e_2-$only
argument is not strong enough in this setting,
so we use $e_3$ as well as $e_2$ and we exploit Mom-2 Technology.  
In particular, we assume
our hyperbolic manifolds contain no geometric Mom-2
involving ${\mathcal O}(1)$ and ${\mathcal O}(2)$,
and then use $e_2$ and $e_3$ to obtain a lower bound
for the area of the maximal cusp torus.
To use the Mom approach to area, we need to be able to apply Theorem
\ref{thm:mom2mflds} and so we require that $e_2 < 1.5152$.
This is certainly the case here.
From Theorem \ref{thm:areanomom2}, we get a  lower bound for the area of the maximal cusp torus of 
\begin{align*}
& 2 \pi (e_3/2)^2 + 2 \pi \left ({e_3 \over e_2} - {e_3 \over 2} \right)^2
 -{\rm overlap}(e_3/2,e_3/2,e_2) \\
& \qquad - 2 \ {\rm overlap}\left ({e_3 \over 2},{e_3 \over e_2} - {e_3 \over 2}, {1 \over e_2} \right).
\end{align*}
 
We know that increasing $e_3$ while holding the other parameters fixed results in an increase 
in the area bound.   The reason is that increasing $e_3$ while holding the other
 parameters fixed increases the radii of the relevant disks but does
 not affect the distance between the centers of the overlapping disks.
 Now, note that for disks of radius $a,\ b$ with centers $c$ apart, if 
 $a$ is increased then the increase in the overlap is a subset of the annulus
 which constitutes 
 the increase in the size of the $a-$radius disk.  Hence the overall effect
 is an increase in area even after overlap is accounted for.  Now, when using
 Tool 1 we see that it works better as area increases.
 Hence, if Tool 1 works for a particular value of $e_3$ then it works for all larger
 values  because an increase in $e_3$ here leads to an increase in area
 (we are assuming, of course, that $1 \le e_2 \le 1.4$).
In {\tt slopes2}, we fix $e_3$ to be $51/32$ and use Tool 1 to verify the
theorems in this case. Thus, in later routines, we may restrict to
$e_2 \leq e_3 \leq 51/32$.

In {\tt slopes3}, we use Tools 1, 2 and 3 to restrict $e_3$ further:
we eliminate parameter points with $1.5 \leq e_3 \leq 51/32$. Thus, in
later routines, we may assume that $e_2 \leq e_3 \leq 1.5$.
 Tool 2 and Tool 3 are a bit trickier here, in that if Tool 2 or 3 works for a particular
 value of $e_3$ then it is possible that Tool 2 or 3 does not work for a larger
 value of $e_3$.  In fact, increasing area makes it easier for Tool 2 and Tool 3
 to fail.  However,  our program is set up so that if Tool 1 fails then the full relevant range
 of parameter values for $m$ and $h$ are analyzed---in particular, 
the area of the fundamental parallelograms $mh$ will exceed the (lower bound
on the) area given by the above calculation.

We finish up by controlling $e_4$

We are now working with $1 \le e_2 \le 1.4$ and $e_2 \le e_3 \le 1.5$.
The above area estimates are inadequate for our methods of eliminating
parameter points to work on these values.  We need to utilize $e_4$
and we now restrict the $e_4$ parameter values. Thus, we use
Theorem \ref{thm:areanomom3}. It follows by
roughly the same reasoning as in the $e_3$ case that if Theorem \ref{thm:areanomom3}
can be used to eliminate a parameter point, then all larger values of $e_4$  are eliminated too.
In fact, the program {\tt slopes4} eliminates points with $e_4 = 51/32$.
Thus we can now restrict our parameter space analysis to $e_3 \le e_4 \le 51/32$.

We have reduced to the following compact parameter space:

$1 \le e_2 \le 1.4$

$e_2 \le e_3 \le 1.5$

$e_3 \le e_4 \le 51/32$

$e_2 \le m \le 2.5$

$1.7 \le h \le 4.0$

$0.0 \le t \le 0.5$

The rigorous analysis of this parameter space using Tools 1,2, and 3 is straightforward.
The relevant program is {\tt slopes5}.

We now discuss computational issues and responses arising from our parameter space analysis.
The computer code was written in C++.
We use interval arithmetic in the form of {\tt Doubles}.  That is, 
we replace {\tt doubles} by intervals and then develop an
arithmetic for these intervals.  Specifically, a {\tt Double}
$x$ is an interval described by a pair of {\tt doubles}
$x = (x.value, x.error)$ where the center of the interval is $x.value$
and the radius of the interval is $x.error$.  Following \cite{GMT1}, we
construct an arithmetic for {\tt Doubles}  whereby if two numbers are
contained within a couple of {\tt Doubles} $x,\ y$ then, for example, the 
{\tt Double} which is the sum
of $x$ and $y$ will contain the actual sum of the original two numbers.

In \cite{GMT1} {\it AffApprox's} are used because of the need for speed.
Our needs are considerably less, and we decided to sneak by using
{\tt Doubles}.  However, because we have 6 parameters, the time 
constraints were still significant, and we resorted to some tricks to
control the time constraints. 

The outer loop of the program {\tt slopes5} corresponds to the  parameter $e_2$.  Next comes the loop
corresponding to $e_3$, then $e_4$.  At this point, we compute a lower bound for 
area using Theorems \ref{thm:areanomom3} and \ref{thm:bonusball}.  The fourth loop corresponds to 
$m$, the minimum translation length.  Given the {\tt Double} $m$ we can use our area bound to compute a lower bound
for $h$.  At this point, we use the functions {\tt CrudeSlopeBound} and
{\tt CrudeIntBound} which determine an upper bound for the number of
exceptional slopes and an upper bound for the maximal intersection number between
exceptional slopes, where both bounds are independent of the parameter $t$.  
If we can eliminate a sub-box of parameter points for
all values of $t$ by using {\tt CrudeSlopeBound} and {\tt CrudeIntBound}, which is a version of Tool 1, then 
we have eliminated the associated sub-boxes (same $e_2, e_3, e_4, m$) 
with all larger values of $h$ as well. When successful,
this is lightning fast, because it is working on a 4-dimensional parameter space.

When {\tt CrudeSlopeBound} and {\tt CrudeIntBound} do not eliminate a sub-box of parameter points then
we introduce the parameter $t$ and do a precise count of the possible number
of exceptional slopes for parameter points for the sub-box in question by using
the functions {\tt FancySlopeBound} and {\tt FancyIntBound}, which are also versions of 
Tool 1.  Again, if we eliminate  a sub-box by this approach, then we have also eliminated the associated sub-boxes
with larger $h$ as well.

If {\tt FancySlopeBound} and {\tt FancyIntBound} do not eliminate a sub-box, then we want to use Tools 
2 and 3.  However, these tools don't automatically eliminate larger $h$ values.
Thus, we turn $h$ into the sixth parameter and use the 3 Tools to
eliminate sub-boxes.

Another technique we used to gain speed was to note that $e_4$ is only
used to get a lower bound on the area of the cusp torus.  Hence,
when a particular {\tt Double} $e_4$ produced enough area (in conjunction with 
the $e_2$ and the $e_3$) then if the next $e_4$ value produced a lower bound 
for area which was at least as large, then we could simply move on to the next
$e_4$ parameter value without bothering with the Tool analysis.  The subtle point here
is that because we are only interested in lower bounds here, we compared the low value
of the intervals in question.

\section{The manifolds in Figure 1}

In this section, we explain the method we used to prove the
following result.

\begin{theorem}
Let $M$ be a compact orientable
hyperbolic 3-manifold with boundary a torus, and which is obtained by Dehn
filling one of the manifolds in Figure \ref{fig:mommflds}. Then the number
of exceptional slopes on $M$ is at most $7$, provided
$M$ is not homeomorphic to $m003$, $m004$ or $m009$.
Moreover, the intersection number between two exceptional slopes
is at most 5, provided $M$ is not homeomorphic to $m003$, $m004$, $m006$ or $m009$.
\label{thm:momcases}
\end{theorem}

To prove this, we could,
in principle, apply the following result to each
of the manifolds in Figure \ref{fig:mommflds}.

\begin{theorem}
Let $N$ be a compact orientable 3-manifold, the interior of which admits
a finite-volume hyperbolic structure. Then 
there is an algorithm to determine all collections of
slopes $(s_0, \dots, s_n)$, with one $s_i$ on each component
of $\partial N$, such that $N(s_0, \dots, s_n)$ is not
hyperbolic.
\label{thm:algorithm}
\end{theorem}

We will not include a proof of this, since it is
a fairly standard application of known algorithms.
It relies on the Casson-Manning
algorithm \cite{Man1} for finding a hyperbolic structure on the interior of
a compact orientable 3-manifold, if one exists. It also uses
normal surface theory algorithms, including
a method for computing the JSJ decomposition of
a manifold \cite{JT1} and the 3-sphere recognition
algorithm of Rubinstein and Thompson \cite{Tho1}. Thus,
it is far from practical. 

Instead, we developed a practical
procedure which creates a set $E$ of slopes $(s_0, \dots, s_n)$ that contains
all the exceptional surgeries. It may be the case that 
$E$ contains some non-exceptional surgeries, but 
these are probably rather rare. 

So, let $N$ be one of the 3-manifolds in Figure \ref{fig:mommflds}.
We will deal later with the unique 3-cusped manifold in Figure \ref{fig:mommflds},
$s776$. We focus now on the case where $N$ has two boundary components.
However, it is clear that this procedure
could be extended to deal with manifolds with more
boundary components.

The procedure relied on the program Snap \cite{Goo1}, which computes hyperbolic
structures on 3-manifolds. Its {\tt verify} function uses exact arithmetic
based on algebraic numbers.
Thus, if the program finds a verified hyperbolic structure, then this is
indeed the correct one.

We first used Snap to find the hyperbolic structure
on $N$. We then used Snap to determine a maximal horoball
neighbourhood of the cusps, with the property that the neighbourhoods
of the two cusps have equal volumes. If $s_0$ and $s_1$ are slopes
on distinct boundary components of $N$ that both have length more
than $6$ with respect to this horoball neighbourhood of the cusps,
then by the 6-theorem and the solution to the geometrisation conjecture,
$N(s_0,s_1)$ is hyperbolic. Thus, we used Snap to determine all the slopes
on the cusp tori with length at most 6.1. (We used 6.1, rather than 6, in order not
to worry about slopes with length precisely 6.) For each such
slope $s_0$ or $s_1$, we needed to determine whether or
not $N(s_0, -)$ (or $N(-,s_1)$) has a hyperbolic structure. (Here, 
$N(s_0,-)$ denotes the manifold obtained by Dehn filling $N$ along
the slope $s_0$, but leaving the second boundary torus unfilled.)
If $N(s_0,-)$ is not hyperbolic, then we declare that $\{ s_0 \} \times S_1$
lies in $E$, where $S_1$ is the set of all slopes on
the other component of $\partial M$. (Recall that $E$ is the set of
slopes that we are aiming to construct, which contains all the exceptional surgeries.)
Similarly, if $N(-,s_1)$ is not hyperbolic, we declare that $S_0 \times \{ s_1 \}$
lies in $E$, where $S_0$ is the set of all slopes on the first
boundary torus. The practical
method of demonstrating that $N(s_0, -)$ or $N(-, s_1)$ was not
hyperbolic was somewhat {\sl ad hoc}, and is described in
more detail below.

If $N(s_0, -)$ is hyperbolic, then we used Snap to find this hyperbolic
structure and to determine a maximal horoball
neighbourhood of its cusp. We then used Snap to find all
slopes with length at most 6.1 on this neighbourhood.
If $s_1$ is a slope on this cusp with length more
than $6.1$, then $N(s_0, s_1)$ is hyperbolic, by the 6-theorem and the
solution to the geometrisation conjecture. We then
considered the slopes $s_1$ with length less than $6.1$,
and used Snap to search for a hyperbolic structure
on $N(s_0, s_1)$. If it could not find one, we
included $(s_0, s_1)$ in $E$. 

We then performed a similar procedure for the hyperbolic manifolds
$N(-, s_1)$, but with the roles of the first and second boundary
components swapped.

Occasionally, slightly indirect methods were required to establish that certain
slopes were non-exceptional. An example is the surgery $s785((4,1),(3,1))$.
Snap asserts that this manifold is hyperbolic. But unfortunately, it cannot
verify this using exact arithmetic. However, it can show that the manifold $s785(-,(3,1))$
is isometric to $m222$ and that the isometry preserves Snap's
co-ordinate systems for the boundary tori. It can also show that
the length of the slope $(4,1)$ on $m222$ is more than $6.1$. Hence,
$s785((4,1),(3,1))$ admits a hyperbolic structure.

In principle, some surgeries in this set $E$ may not be exceptional.
This can happen in two ways. The first situation
arises when $N(s_0, -)$ is non-hyperbolic.
For example, $N(s_0, -)$ may have non-trivial
JSJ decomposition, but $N(s_0, s_1)$ may still be
hyperbolic for some $s_1$. Theorem \ref{thm:algorithm} provides a theoretical algorithm
for finding all such $s_1$, but this is not practical. However, the set of such
$s_1$ will probably be rather sparse, and so it did not seem too wasteful
to include them in the set $E$. The second way
that $E$ may be too large arose in the search
for a hyperbolic structure on $N(s_0, s_1)$,
as a Dehn filling of a fixed hyperbolic $N(s_0, -)$. 
It quite often happens that Snap does not find a
hyperbolic structure when one exists. In practice,
one may need to use Snap to retriangulate the
manifold several times. Even then, manifolds
where Snap fails to find a hyperbolic structure
may yet have one. Thus, $E$ may be somewhat larger than the
actual set of exceptional surgeries. However, it
is certainly the case that all exceptional surgeries
lie in $E$. And the procedure was discerning enough
to produce Theorem \ref{thm:momcases}.

It remains to describe how we dealt with the cases
where $N(s_0, -)$ (or $N(-,s_1)$) was not hyperbolic. Here, it is
important to prove that $N(s_0, -)$ definitely does not admit a hyperbolic
structure, rather than simply declaring that Snap could
not find one. This is because, otherwise,
$N(s_0, -)$ might in fact have had a hyperbolic structure,
and then be a potential counter-example to Theorem \ref{thm:momcases}.

In practice, we proved that $N(s_0, -)$ was not
hyperbolic using a variety of methods, all of
which utilised a presentation for $\pi_1(N(s_0, -))$.
Snap can produce not only this presentation,
but also give generators for the peripheral
subgroups. We used the following four techniques for proving
that $N(s_0,-)$ is not hyperbolic.
\begin{enumerate}
\item If this group presentation is obviously that
of a non-trivial free product, then $N(s_0,-)$ is
not hyperbolic.
\item If the existence of a non-trivial
centre is obvious from this presentation, then
again $N(s_0, -)$ is not hyperbolic. In each case,
an element of the group was found which commuted with every
generator, and which could be seen to be homologically
non-trivial.
\item If the group contains a commuting pair of
elements that lie in neither a cyclic subgroup nor a peripheral
subgroup, the manifold is not hyperbolic.
In practice, both of the required properties
of these two elements were verified homologically.
\item If the presentation has two generators
$a$ and $b$, and the same non-trivial power of
$a$ appears in all the relations, the manifold is
not hyperbolic. We now supply a proof of this.
\end{enumerate}

\noindent {\sl Proof.} Suppose the presentation
is $\langle a, b | w_1(a^k, b), \dots, w_r(a^k,b) \rangle$,
where $k > 1$ and the $w_i(a^k,b)$ are words in $a^k$ and $b$.
This group is then an amalgamated free product
$$\langle a \rangle \ast_{a^k = c} 
\langle c, b | w_1(c, b), \dots, w_r(c,b) \rangle.$$
Suppose first that this is a trivial amalgamated
free product. Then the amalgamating subgroup
must be the whole of one of the factors.
It cannot be the first factor, since 
$a^k$ is a proper subgroup of $\langle a \rangle$.
Thus, it must be the second factor, which implies
that the group is $\langle a \rangle$, which is
cyclic. Hence, in this case, the manifold is not
hyperbolic. On the other hand, if this is a non-trivial
amalgamated free product, then the manifold has non-trivial
JSJ decomposition, because the amalgamating subgroup
is cyclic. Again, this implies that the manifold
is not hyperbolic. \qed\

The following table gives a summary of where
these methods were applied. In the first column,
the manifold $N$ is given. In the second column,
Snap's label for the relevant boundary component is given, 
which is either $0$ or $1$. In the third column,
the slope $s_0$ or $s_1$ is shown, in the co-ordinates given by Snap.
The final column gives a number between 1 and 4, according
to the method used to prove non-hyperbolicity.
It turns out that, in all the manifolds we considered,
except $s785$, there exists an isometry of the manifold, swapping the
cusps, and preserving Snap's co-ordinate system for
the slopes. Thus, in all the cases except $s785$, we only
give the slopes on the boundary component labelled $0$.

We briefly mention $s780((1,1),-)$ since this was a slightly tricky case.
Snap provides the following presentation of its fundamental group:
$$\langle a,b \mid ab^{-1}a^{-1}ba^2ba^{-1}b^{-1}ab^{2} \rangle.$$
None of the above four methods can be obviously applied here. However,
after re-triangulating, we obtain a new presentation
$$\langle a,b,c \mid bc^{-1}bc^2, ab a^{-1}b \rangle.$$
The second relation is that of the Klein bottle fundamental group.
Thus, $a^2$ and $b$ commute. They generate a subgroup of $H_1(M)$
isomorphic to ${\Bbb Z} \oplus {\Bbb Z}/2$, which is not cyclic.
The peripheral subgroup is $\langle bc, ca^{-1}cb^{-1}ca^{-1}c^{-1} \rangle$,
which, in $H_1(M)$ is exactly the subgroup generated by $a^2$ and $b$.
So, this does not immediately lead to a contradiction. However, if $a^2$
were peripheral, then so would $a$ be. But $a$ does not lie in the image of
$H_1(\partial M)$ in $H_1(M)$.

\begin{center}
\begin{tabular}{|c|c|c|c||c|c|c|c|}
\hline\hline
$m412$ & 0 & (1,0)  & 1 &  $s785$ & 0 & (1,0)  & 2 \\
     & 0 & (-1,1) & 3 &       & 0 & (-1,1) & 4 \\
     & 0 & (0,1)  & 2 &       & 0 & (0,1)  & 2 \\
     & 0 & (1,1)  & 3 &       & 0 & (1,1)  & 4 \\
$s596$ & 0 & (1,0)  & 1 &       & 1 & (1,0)  & 2 \\
     & 0 & (-1,1) & 2 &       & 1 & (-1,1) & 2 \\
     & 0 & (0,1)  & 2 &       & 1 & (0,1)  & 2 \\
     & 0 & (-1,2) & 4 &       & 1 & (1,1)  & 2 \\
$s647$ & 0 & (1,0)  & 2 &  $s898$ & 0 & (1,0)  & 2 \\
     & 0 & (-1,1) & 2 &       & 0 & (-1,1) & 4 \\
     & 0 & (0,1)  & 2 &       & 0 & (0,1)  & 2 \\
     & 0 & (1,1)  & 4 &       & 0 & (1,1)  & 2 \\
$s774$ & 0 & (1,0)  & 1 &  $s959$ & 0 & (1,0)  & 2 \\
     & 0 & (-1,1) & 4 &       & 0 & (-1,1) & 2 \\
     & 0 & (0,1)  & 2 &       & 0 & (0,1)  & 2 \\
$s780$ & 0 & (1,0)  & 2 &	      &  &         & \\
     & 0 & (1,1)  & 3 &       &  &         & \\ 
\hline\hline
\end{tabular}
\end{center}

We applied the above procedure to the manifolds
$m412$, $s596$, $s647$, $s774$, $s780$, $s785$, $s898$ and $s959$. 
We were able to construct a set $E$ for each manifold $N$,
and thereby verify Theorem \ref{thm:momcases} for any hyperbolic manifold
$M$ obtained by Dehn filling $N$.

As an example, we include here the set $E$ for $m412$.
In this case, $E$ is
\begin{align*}
&(1,0) \times S_1, \qquad
(-1,1) \times S_1, \qquad
(0,1) \times S_1, \qquad 
(1,1) \times S_1, \\
&S_0 \times (1,0), \qquad
S_0 \times (-1,1), \qquad
S_0 \times (0,1), \qquad
S_0 \times (1,1),
\end{align*}
together with the pairs $(s_0, s_1)$ in the table below marked with an x.

\begin{center}
\begin{tabular}{|c|c|c|c|c|c|c|c|c|}
\hline\hline
&  		 (-2,1) & (2,1) & (3,1) & (4,1) & (5,1) & (-1,2) & (-1,3) \\
\hline
(-2,1) &	 	&	&	&	&	& x	& x \\
\hline
(2,1) & 		& x 	& x	& x	& x 	&	& \\
\hline
(3,1) &			& x	& x 	&	&	&	& \\
\hline
(4,1) & 		& x 	&	&	&	&	& \\
\hline
(5,1) & 		& x 	&	&	&	&	& \\
\hline
(-1,2) &  	 x 	&	&	&	&	& 	& \\
\hline
(-1,3) &	 x 	&	&	&	&	&	& \\
\hline\hline
\end{tabular}
\end{center}

Here, the slopes along the top are on boundary component $0$,
and those down the left are on boundary component $1$.
So, for instance, $M((2,1),-)$ is a hyperbolic manifold with at most $8$ exceptional surgeries:
$(1,0)$, $(-1,1)$, $(0,1)$, $(1,1)$, $(2,1)$, $(3,1)$, $(4,1)$, $(5,1)$.
In fact, $M((2,1),-)$ is $m009$, which is known to have precisely 8 exceptional surgeries.

The sets $E$ for the manifolds  $s596$, $s647$, $s774$, $s780$, $s785$, $s898$ and $s959$
can be found in the appendices. 

For each manifold with more than $7$ exceptional slopes, or where the
intersection number between two exceptional slopes is more than $5$, we need
to prove that it is one of the manifolds in Theorem \ref{thm:momcases}.
This is achieved using Snap's {\tt identify} command.

This leaves one remaining manifold from Figure \ref{fig:mommflds}, $s776$. 
This has three boundary components. So, although the above procedure can,
in principle, be applied here, it is less practical.
Fortunately, Martelli and Petronio \cite{MP1} have determined
all the exceptional surgeries on this manifold $M$.
In particular, Corollary A.6 in \cite{MP1} implies the first part of Theorem \ref{thm:momcases}
in this case. Also, by going through the explicit list of
exceptional surgeries in \cite{MP1}, it is possible to verify the
second part of Theorem \ref{thm:momcases}. \qed\

Note that the manifolds $m003$, $m004$, $m006$ and $m009$ referred to in Theorem
\ref{thm:momcases} all satisfy Theorems \ref{thm:main1} and \ref{thm:main2}.
This therefore completes the proof of these theorems. \qed\

\appendix
\section{Surgeries on $s596$}

\begin{align*}
&(1,0) \times S_1, \qquad
(-1,1) \times S_1, \qquad
(0,1) \times S_1, \qquad 
(-1,2) \times S_1, \\
&S_0 \times (1,0), \qquad
S_0 \times (-1,1), \qquad
S_0 \times (0,1), \qquad
S_0 \times (-1,2),
\end{align*}

\begin{center}
\begin{tabular}{|c|c|c|c|c|}
\hline\hline
&  		 (-2,1) & (1,1) & (-2,3) & (-1,3) \\
\hline
(-2,1) &	x	&	& x	&	\\
\hline
(1,1) & 		& x 	& 	& x	\\
\hline
(-2,3) &	x	&	& 	&  	\\
\hline
(-1,3) & 		& x 	&	&	\\
\hline\hline
\end{tabular}
\end{center}

\section{Surgeries on $s647$}

\begin{align*}
&(1,0) \times S_1, \qquad
(-1,1) \times S_1, \qquad
(0,1) \times S_1, \qquad 
(1,1) \times S_1, \\
&S_0 \times (1,0), \qquad
S_0 \times (-1,1), \qquad
S_0 \times (0,1), \qquad
S_0 \times (1,1),
\end{align*}

\begin{center}
\begin{tabular}{|c|c|c|c|c|c|c|c|c|c|c|}
\hline\hline
&  		 (-7,1) & (-6,1) & (-5,1) & (-4,1) & (-3,1) & (-2,1) & (2,1) & (-1,2) & (1,2) & (-2,3)\\
\hline
(-7,1) &		&	& 	&	&	& x	&  		&	&	&	\\
\hline
(-6,1) & 		&  	& 	& 	&	& x	& 		&	&	&	\\
\hline
(-5,1) &		&	& 	&  	&	& x	& 		&	&	&	\\
\hline
(-4,1) & 		&  	&	&	& 	& x	& 		&	&	&	\\
\hline
(-3,1) & 		&  	&	&	& x	& x	& 		&	&	&	\\
\hline
(-2,1) & 	x	& x  	& x	& x	& x	& x	& 		&	&	&	\\
\hline
(2,1) & 		&   	& 	& 	& 	& 	& 		& x	&	& x	\\
\hline
(-1,2) & 		&  	&	&	& 	& 	& x		&	& x	&	\\
\hline
(1,2) & 		&  	&	&	& 	& 	& 		& x	& x	&	\\
\hline
(-2,3) & 		&  	&	&	& 	& 	& x		&	&	&	\\
\hline\hline
\end{tabular}
\end{center}

\section{Surgeries on $s774$}

\begin{align*}
&(1,0) \times S_1, \qquad
(-1,1) \times S_1, \qquad
(0,1) \times S_1, \\
&S_0 \times (1,0), \qquad
S_0 \times (-1,1), \qquad
S_0 \times (0,1), 
\end{align*}

\begin{center}
\begin{tabular}{|c|c|c|c|c|c|c|c|c|c|}
\hline\hline
&  		 (-5,1) & (-4,1) & (-3,1) & (-2,1) & (1,1) & (2,1) & (3,1) & (-1,2) & (-1,3) \\
\hline
(-5,1) &		&	& 	& x	&	& 	&  		&	&	\\
\hline
(-4,1) & 		&  	& 	& x	&	& 	& 		&	&	\\
\hline
(-3,1) &		&	& x	& x 	&	& 	& 		&	&	\\
\hline
(-2,1) & 	x	& x  	& x	& x	& 	& 	& 		&	&	\\
\hline
(1,1) & 		&  	&	&	& x	& 	& 		& x	& x	\\
\hline
(2,1) & 		&   	& 	& 	& 	& 	& 		& x	&	\\
\hline
(3,1) & 		&  	&	&	& 	& 	& 		& x	&	\\
\hline
(-1,2) & 		&  	&	&	& x	& x	& x		&	& 	\\
\hline
(-1,3) & 		&  	&	&	& x 	& 	& 		&	&	\\
\hline\hline
\end{tabular}
\end{center}

\section{Surgeries on $s780$}

\begin{align*}
&(1,0) \times S_1, \qquad
(1,1) \times S_1, \\
&S_0 \times (1,0), \qquad
S_0 \times (1,1), 
\end{align*}

\begin{center}
\begin{tabular}{|c|c|c|c|c|c|c|c|c|c|c|c|c|}
\hline\hline
&  		 (-7,1) & (-6,1) & (-5,1) & (-4,1) & (-3,1) & (-2,1) & (-1,1) & (0,1) 	& (2,1) & (3,1) & (4,1) & (5,1)\\
\hline
(-7,1) &		&	& 	&	&	& 	&  		& x	&	&	&	&	\\
\hline
(-6,1) & 		&  	& 	& 	&	& 	& 		& x	&	&	&	&	\\
\hline
(-5,1) &		&	& 	&  	&	& 	& x		& x	&	&	&	&	\\
\hline
(-4,1) & 		&  	&	&	& 	& x	& x		& x	&	&	&	&	\\
\hline
(-3,1) & 		&  	&	&	& x	& x	& x		& x	&	&	&	&	\\
\hline
(-2,1) & 		&   	& 	& x	& x	& x 	& x 		& x	&	&	&	&	\\
\hline
(-1,1) & 		&   	& x 	& x	& x	& x	& x		& x 	&	& 	&	&	\\
\hline
(0,1) & 	x	& x 	& x	& x 	& x	& x 	& x		& x	&	&	&	&	\\
\hline
(2,1) & 		&  	&	&	& 	& 	& 		&	& x	& x	& x	& x	\\
\hline
(3,1) & 		&  	&	&	& 	& 	& 		&	& x	& x	&	&	\\
\hline
(4,1) & 		&  	&	&	& 	& 	& 		&	& x	&	&	&	\\
\hline
(5,1) & 		&  	&	&	& 	& 	& 		&	& x	&	&	&	\\
\hline\hline
\end{tabular}
\end{center}

\section{Surgeries on $s785$}

\begin{align*}
&(1,0) \times S_1, \qquad
(-1,1) \times S_1, \qquad
(0,1) \times S_1, \qquad
(1,1) \times S_1, \\
&S_0 \times (1,0), \qquad
S_0 \times (-1,1), \qquad
S_0 \times (0,1), \qquad
S_0 \times (1,1), 
\end{align*}

\begin{center}
\begin{tabular}{|c|c|c|c|c|c|}
\hline\hline
&  		 (-2,1) & (2,1) & (3,1) & (-1,2) & (1,2) \\
\hline
(-2,1) &		& x	& x	&	& x	\\
\hline
(2,1) & 		& x 	& x 	& 	&	\\
\hline
(-1,2) &		& 	&  	& x 	&	\\
\hline
(1,2) & x 		&  	&	& 	&	\\
\hline
(2,3) & x		&  	&	&	&	\\
\hline\hline
\end{tabular}
\end{center}

\section{Surgeries on $s898$}

\begin{align*}
&(1,0) \times S_1, \qquad
(-1,1) \times S_1, \qquad
(0,1) \times S_1, \qquad
(1,1) \times S_1, \\
&S_0 \times (1,0), \qquad
S_0 \times (-1,1), \qquad
S_0 \times (0,1), \qquad
S_0 \times (1,1), 
\end{align*}

\begin{center}
\begin{tabular}{|c|c|c|c|c|c|}
\hline\hline
  	& (-3,1)& (-2,1) & (2,1) & (-1,2) & (1,2) \\
\hline
(-3,1) & 	& x	&	&	&	\\
\hline
(-2,1) & x	& x	& x	& 	&	\\
\hline
(2,1) &  	& x	& x  	& x 	& 	\\
\hline
(-1,2) &	&	& x	&  	&  	\\
\hline
(1,2) &  	&	&  	&	& x 	\\
\hline\hline
\end{tabular}
\end{center}

\section{Surgeries on $s959$}

\begin{align*}
&(1,0) \times S_1, \qquad
(-1,1) \times S_1, \qquad
(0,1) \times S_1, \\
&S_0 \times (1,0), \qquad
S_0 \times (-1,1), \qquad
S_0 \times (0,1) 
\end{align*}

\begin{center}
\begin{tabular}{|c|c|c|c|c|c|c|}
\hline\hline
&  		 (-3,1) & (-2,1) & (1,1) & (-1,2) & (1,2) & (-2,3) \\
\hline
(-3,1) &		& 	& 	& x	& 	&	\\
\hline
(-2,1) & 		& x 	&  	& x	&	& x	\\
\hline
(1,1) &			& 	& x 	& x 	& x	&	\\
\hline
(-1,2) &  	x	& x 	& x	& 	& 	&	\\
\hline
(1,2) & 		&  	& x 	&	&	&	\\
\hline
(-2,3) & 		& x 	&	&	&	&	\\
\hline\hline
\end{tabular}
\end{center}


\begin{thebibliography}{99}
\bibitem[28]{Per1} G. Perelman, {\sl The entropy formula for the Ricci flow and its geometric applications},
arxiv:math.DG/0211159
\end{thebibliography}
\end{document}